\documentclass[12pt]{elsarticle}

\usepackage{lineno,hyperref}
\modulolinenumbers[5]
\usepackage{enumitem}
\usepackage{array}
\usepackage{tabularx} 
\usepackage{adjustbox} 
\usepackage{makecell}
\usepackage{amsmath}
\usepackage{setspace}
\usepackage{float}
\usepackage{graphicx}
\usepackage{subfigure}
\usepackage{geometry}
\usepackage{amsmath}
\usepackage{amsfonts}
\usepackage{amssymb}

\journal{Journal of \LaTeX\ Templates}

\begin{document}

\begin{frontmatter}

\title{\Large \textbf{Solving linear objective optimization problem subjected to novel max-min fuzzy relational equalities as a generalization of the vertex cover problem}}

\author[1]{Amin Ghodousian\corref{mycorrespondingauthor}}
\cortext[mycorrespondingauthor]{Corresponding author}
\ead{a.ghodousian@ut.ac.ir}
\author[2]{Mahdi Mollakazemiha}

\address[1]{School of Engineering Science, College of Engineering, University of Tehran, Tehran, Iran.}

\address[2]{Faculty of Mathematical Science, Shahid Beheshti University , Tehran, Iran.}

\begin{large}
\begin{abstract}
\setstretch{1.15}
This paper considers the linear objective function optimization with respect to a novel system of fuzzy relation equations, where the fuzzy compositions are defined by the minimum t-norm. It is proved that the feasible solution set is formed as a union of the finite number of closed convex cells. Some necessary and sufficient conditions are presented to conceptualize the feasibility of the problem. Moreover, seven rules are introduced with the aim of simplifying the original problem, and then an algorithm is accordingly presented to find a global optimum. It is shown that the original problem in a special case is reduced to the well-known minimum vertex cover problem. Finally, an example is described to illustrate the proposed algorithm.

\end{abstract}
\end{large}

\end{frontmatter}

\begin{large}

\section{Introduction}
\onehalfspacing

\noindent
The theory of fuzzy relational equations (FRE) as a generalized version of Boolean
relation equations was firstly proposed by Sanchez and was applied to problems
related to the medical diagnosis [40]. Pedrycz categorized and extended two ways
of the generalizations of FRE in terms of sets under discussion and various
operations which are taken account [36]. Since then, FRE was applied in many other
fields such as fuzzy control, prediction of fuzzy systems, fuzzy decision making,
fuzzy pattern recognition, image compression and reconstruction, fuzzy clustering
and so on. Generally, when rules of inference are applied and their corresponding consequences are known, the problem of determining antecedents is simplified
and mathematically reduced to solving an FRE [34]. Nowadays, it is well known that
many of the issues associated to the body knowledge can be treated as FRE
problems [35]. Because of important applications in various practical fields, many
scholars have focused on the theoretical research of FRE, including its resolution
approach and some specific optimization problems with FRE constraints. \\
The solvability identification and finding set of solutions are the primary, and the
most fundamental, matter concerning the FRE problems. Di Nola et al. proved that
the solution set of FRE (if it is nonempty), defined by continuous max-t-norm
composition is often a non-convex set. This non-convex set is completely
determined by one maximum solution and a finite number of minimal solutions [6].
Such non-convexity property is one of two bottlenecks making a major contribution
towards an increase in complexity of FRE-related problems, particularly, in the
optimization problems subjected to a system of fuzzy relations. Another bottleneck
point is concerned with detecting the minimal solutions for FREs. Chen and Wang
[2] presented an algorithm for obtaining the logical representation of all minimal
solutions and deduced that a polynomial-time algorithm with the ability to find all
minimal solutions of FRE (with max-min composition) may not exist. Also,
Markovskii showed that solving max-product FRE is closely related to the covering
problem which is a type of NP-hard problem [33]. \\
In fact, the same result holds
true for a more general t-norms instead of the minimum and product operators
[3,28,29]. Over the past decades, the solvability of FRE which is defined using
different max-t compositions have been investigated by many researchers
[16,17,19,38,42,45,46,49,51]. Moreover, some other researchers have worked on
introducing novel concept and at times improving some of the existing theoretical
aspects and applications of fuzzy relational inequalities (FRI) [14,18-20,25,53]. For example, Li and Yang [25] studied FRI with addition-min composition and presented an algorithm to search for minimal solutions. They applied FRI to data transmission mechanism in a BitTorrent-like Peer-to-Peer file sharing systems. In [14], the authors focused on the study of a mixed fuzzy system formed by two fuzzy
relational inequalities $A \phi x \le b^{1}$ and $D \phi x \ge b^{2}$, where $\phi$ is an operator with (closed) convex solutions. \\
The optimization problem subject to FRE and FRI is one of the most interesting and
on-going research topics amongst similar problems [1,10,13,14,16,22,24,30,39
,43,47,53]. Many methods were designed based on the translation of the main problem into an integer linear programming problem which is then solved using well-developed techniques. On the contrary, other algorithms benefit the resolution of the feasible region, some necessary and sufficient conditions for the optimality and simplification processes. The most methods of this category are based on analytical results provided mainly by Sanchez [41] and Pedrycz [37]. For instance, Fang and Li converted a linear optimization problem subjected to FRE constraints with max-min operation into an integer programming problem and solved it by a branch-and-bound method using jump-tracking technique [11]. Wu et al. worked on improvement of the method employed by Fang and Li; this was done by decreasing the search domain and presented a simplification process by
three rules which resulted from a necessary condition [48]. Chang and Shieh
presented new theoretical results concerning the linear optimization problem
constrained by fuzzy max–min relation equations [1]. They improved an upper
bound on the optimal objective value, some rules for simplifying the problem and
proposed a rule for reducing the solution tree. In [23], an application of optimizing
the linear objective with max-min composition was employed for the streaming
media provider seeking a minimum cost while fulfilling the requirements assumed
by a three-tier framework. Linear optimization problem was further investigated
by numerous scholars focusing on max-product operation [22,32]. Loetamonphong
and Fang defined two sub-problems by separating negative from non-negative
coefficients in the objective function, and then obtained an optimal solution by
combining the optimal solutions of the two sub-problems [32]. Moreover,
generalizations of the linear optimization problem with respect to FRE have been
studied; this was done through replacement of max-min and max-product
compositions with different fuzzy compositions such as max-average composition
[47] or max-t-norm composition [16,17,19,21,24,43]. For example, Li and Fang
solved the linear optimization problem subjected to a system of sup-t equations by
reducing it to a 0-1 integer optimization problem [24]. In [21], a method was
presented for solving linear optimization problems with the max-Archimedean tnorm fuzzy relation equation constraint. In [43], the authors solved the same
problem whit continuous Archimedean t-norm, and to obtain some optimal
variables, they used the covering problem rather than the branch-and-bound
methods. \\ 
Recently, many interesting forms of generalizations of the linear programming
applied to the system of fuzzy relations have been introduced, and developed
based on composite operations used in FRE, fuzzy relations used in the definition
of the constraints, some developments on the objective function of the problems
and other ideas [5,7,12,19,27,30,50,54]. For example, Wu et al. represented an
efficient method to optimize a linear fractional programming problem under FRE
with max-Archimedean t-norm composition [50]. Dempe and Ruziyeva generalized
the fuzzy linear optimization problem by considering fuzzy coefficients [5]. In
addition, Dubey et al. studied linear programming problems involving interval
uncertainty modeled using intuitionistic fuzzy set [7]. Yang [54] studied the optimal
solution of minimizing a linear objective function subject to a FRI where the
constraints defined as $\sum_{j=1}^{n}\min\{ a_{ij}, x_{j} \} \ge b_{i}$ for $i = 1,2,...,m$. Also, in [53], the authors
introduced the latticized linear programming problem subject to max-product fuzzy
relation inequalities with application in the optimization management model for
wireless communication emission base stations. The latticized linear programming
problem was defined by minimizing the objective function $z(x) = \min_{j=1}^{n}\{ 
x_{j} \}$ subject to the feasible region $X(A,b) = \{ x \in [0,1]^{n} : A \circ x \ge b \}$ where “$\circ$” denotes fuzzy maxproduct composition. They also presented an algorithm based on the resolution of the feasible region. \\
The concept of FRE was also generalized to a so-called bipolar fuzzy relation
equation. Bipolarity exists widely in human understanding of information and
preference [8]. Dubois and Prade provided an overview of the asymmetric bipolar
representation of positive and negative information in possibility theory [9]. They
shown that the possibility theory framework is convenient for handling bipolar
representations, that was applied to distinguish between negative and positive
information in preference modeling [8,9]. The linear optimization of bipolar FRE
was also the focus of study carried out by some researchers where FRE was defined
with max-min (with application in public awareness of the products for a supplier)
[12,26], max-product [4] and max-Lukasiewicz composition [27,30,52,55]. In [12],
the concept of bipolar FRE was firstly proposed with max-min composition where
the constraints are expressed as $\max_{j=1}^{n} \{ \max\{  \min\{ a_{ij}^{+}, x_{j} \} , \min\{ a_{ij}^{-}, 1 - x_{j} \} \} \}$ for $i=1,2,...,m$, where $a_{ij}^{+} , a_{ij}^{-}, x_{j} \in [0,1]$. Similarly, in [27], the authors introduced a linear
optimization problem subjected to a system of bipolar FRE defined as 
$X(A^{+},A^{-},b) = \{ x \in [0,1]^{m} : x \circ A^{+} \lor \Tilde{x} \circ A^{-} = b \}$, where $\Tilde{x_{i}} = 1 - x_{i} $ for each component of $\Tilde{x} = (\Tilde{x_{i}})_{1 \times m}$ and the notations “$\circ$” and “$\lor$” denote max operation and the max-Lukasiewicz composition, respectively. They translated the original problem into a 0-1 integer linear problem. In a separate, the foregoing bipolar linear optimization problem was solved by an analytical method based on the resolution and some structural properties of the feasible region (using a necessary condition for characterizing an optimal solution and a simplification process for reducing the problem) [30]. However, resolution method for obtaining the complete solution set of bipolar max-Lukasiewicz FRE was not found in the mentioned works [52]. Yang [52] studied bipolar max-Lukasiewicz FRE and showed that the complete solution set of the system is fully determined by finite conservative bipolar paths. \\
The original problem of solving max-min fuzzy relational equations is defined as the following system: \\

\begin{equation}
    \begin{array}{l}
\textbf{A} \circ \textbf{x} = \textbf{b} \\ 
x \in [0,1]^{n}
    \end{array}
\end{equation}

where $\textbf{A} = (a_{ij})_{m \times n}$ and $\textbf{b} = (b_{i})_{m \times 1}$ are the coefficients matrix and the right-handside vector, respectively, whose components belong to the interval [0,1], $"\circ"$ denotes the max-min composition and the constraints mean $\max_{j=1}^{n}\{ \min\{ a_{ij}, x_{j} \} \} = b_{i}$, $i = 1,2,...,m$. Given that \textbf{A} and \textbf{b} are known, the resolution problem is to determine all vectors x such that constraints (1) are satisfied. Figure 1 schematically shows the feasible region of the original max-min FRE problems. 

\begin{figure}[ht]
    \begin{center}
	    \includegraphics[height=7cm]{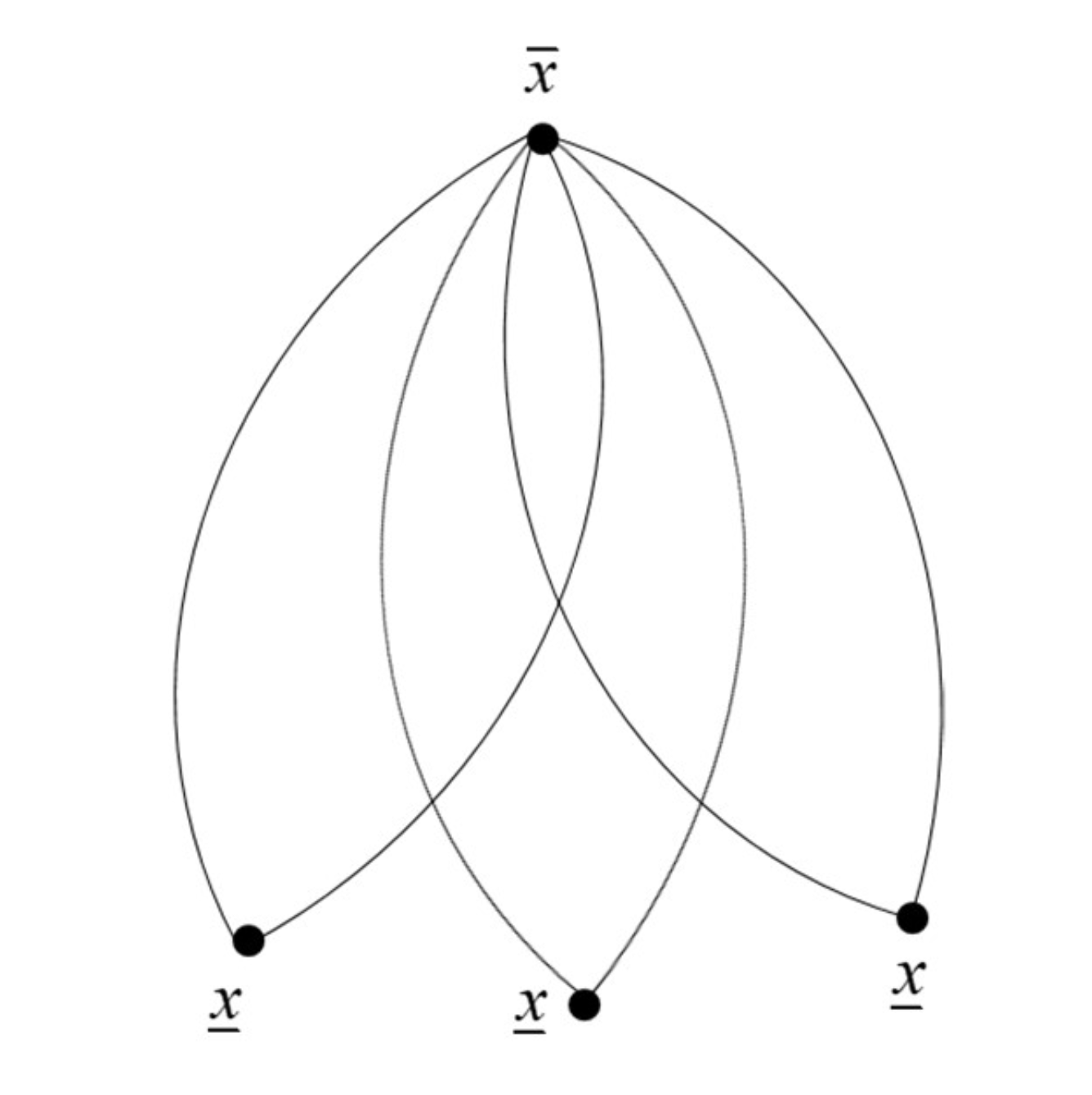}
	    \caption{The domain of original max-min FRE.}
	\end{center}
\end{figure} 

In this paper, we study the following optimization problem:\\

\begin{equation}
    \begin{array}{l}
        \min Z = \textbf{c}^{T} \textbf{x} \\
        \textbf{A} \otimes \textbf{x} = \textbf{b} \\
        \textbf{x} \in [0,1]^{n}
        
    \end{array}
\end{equation}

where $I \in \{1,2,...,m \}$ and $J \in \{1,2,...,n \}$. $\textbf{A} = (a_{ij})_{m \times n}$ is the coefficient matrix such that $a_{ij} \in [0,1]$ $(\forall i \in I$ and $\forall j \in J)$, $\textbf{b} = (b_{i})_{m \times 1}$ is the right-hand-side vector such that $b_{i} \in [0,1]$ $(\forall i \in I)$ and $\textbf{c}$ is a vector in $\mathbf{R}^{n}$. Moreover, if $\textbf{a}_{i}$ denotes the i‘th row of matrix \textbf{A} , then the i‘th constraint of Problem (2) is defined as follows:\\

\begin{equation}
    \begin{array}{l}
        \textbf{a}_{i} \otimes \textbf{x} = \max_{j=1}^{n}\{\min\{ a_{ij}, x_{i}, x_{j}\} \} = b_{i}, ~~~~ , ~~~~ i \in I
    \end{array}
\end{equation}

It should be noted here that without loss of generality, it can be assumed that \textbf{A} is a square matrix (i.e., $m = n$); because otherwise, in the case that $m > n$, $m-n$ additional columns $j_{n+1},...,j_{m}$ with corresponding coefficients $a_{ij_{k}} = 0 ~~ (i \in I$ and $k \in \{n+1,...,m\})$ can be augmented to \textbf{A} , and in the other case that $m < n$, $n-m$ equations $\textbf{a}_{i} \otimes \textbf{x} = \max_{j=1}^{n}\{\min\{ 0, x_{i}, x_{j}\} \} = 0 ~~ (i \in \{ m+1, ..., n  \})$ can be added to the problem. So, hereafter, m is assumed to be equal to n throughout the paper.\\

The rest of the paper is organized as follows. Some preliminary definitions and concepts are presented and the feasible solution set of the problem is subsequently determined in Section 2. Moreover, some necessary and sufficient conditions are derived to identify the feasibility of the problem. In Section 3, seven rules are introduced to accelerate the resolution process by eliminating some special cases from consideration and reducing the size of the problem. Problem (2) is resolved through optimizing of the linear objective function which is considered in Section 4. In addition, the existence of a global optimal solution for the problem is proved in the general case, and it is shown that in the special case where the right-hand- side vector \textbf{b} is the zero vector \textbf{0} and all the coefficients and the variables are assumed to be binary, then all the optimal solutions of the problem are also obtained in binary form. The previously obtained results are summarized as an algorithm and in Section 5, it is shown that the well-known minimum vertex cover problem is actually a special case of the main problem. Finally, in Section 6 an illustration of the general algorithm is provided using an example.\\

\section{Feasible solutions set of the problem}

This section describes the characterization of the feasible region of Problem (2). For this purpose, we firstly determine the feasible solutions set of the i’th constraint of Problem (2), i.e., Relation (3) for each $i \in I$. \\

\textbf{Definition 1.} 
For each $i \in I$, let $S(\textbf{a}_{i},b_{i})$ denotes the feasible solutions set of i’th equation, i.e., $ S(\textbf{a}_{i},b_{i}) = \{ \textbf{x} \in [0,1]^{n}: \textbf{a}_{i} \otimes  \textbf{x} = b_{i}\} = \{ \textbf{x} \in [0,1]^{n} : \max_{j \in J} \{\min\{a_{ij}, x_{i}, x_{j}  \}  \} = b_{i} \}$. Also let $S(\textbf{A}, \textbf{b}) = \{\textbf{x} \in [0,1]^{n}: \textbf{A} \otimes \textbf{x} = \textbf{b} \}$. So, $S(\textbf{A}, \textbf{b})$ denotes the feasible solutions set of Problem (2).\\
A solution $\Bar{\textbf{x}} = (\Bar{x_{1}}, ..., \Bar{x_{n}}) \in S(\textbf{A}, \textbf{b})$ is said to be a maximum solution when $\textbf{x} \le \Bar{\textbf{x}}$ (i.e., $x_{j} \le \Bar{x_{j}}, ~~ \forall j \in J$) for all $\textbf{x} \in S(\textbf{A}, \textbf{b}) $. Similarly, $\underline{\textbf{x}} \in S(\textbf{A}, \textbf{b})  $ is said to be a minimum solution when $\underline{\textbf{x}} \le \textbf{x}$ for all $\textbf{x} \in S(\textbf{A}, \textbf{b}) $. Moreover, a solution $\Bar{\textbf{x}} \in S(\textbf{A}, \textbf{b}) ~ (\underline{\textbf{x}} \in S(\textbf{A}, \textbf{b}))$ is said to be a maximal (minimal) solution when $\Bar{\textbf{x}} \le \textbf{x} ~ (\textbf{x} \le \underline{\textbf{x}})$ implies $\Bar{\textbf{x}} = \textbf{x} ~ (\textbf{x} = \underline{\textbf{x}}) $ for any $\textbf{x} \in S(\textbf{A}, \textbf{b}) $.\\
From Definition 1, it is clear that  $S(\textbf{A}, \textbf{b}) = \bigcap_{i \in I}{S(\textbf{a}_{i},b_{i})}$. Moreover, this definition together with Relation (3) result in Lemma 1 below that provides a necessary and sufficient condition for the feasibility of $S(\textbf{a}_{i},b_{i}), ~ \forall i \in I$.\\

\textbf{Lemma 1.} For a fixed $i \in I$, $\textbf{x} \in S(\textbf{a}_{i},b_{i})$ iff $\textbf{x} \in [0,1]^{n}$ and the following two conditions are satisfied:

\begin{equation}
    \begin{array}{l}
        (i)~~ \min\{a_{ij}, x_{i}, x_{j} \} \le b_{i}, ~~~ \forall j \in J.\\
        (ii)~~ \min\{a_{ij_{0}}, x_{i}, x_{j_{0}} \} = b_{i} ~~~ for ~ some ~ j_{0}\in J.
    \end{array}
\end{equation}

\textbf{Definition 2.} For each $i \in I$, we define $J_{i}^{1}= \{j \in J: a_{ij} > b_{i}  \}$ and $J_{i}^{2}= \{j \in J: a_{ij} = b_{i}  \}$. Also let $J_{i} = J_{i}^{1} \bigcup J_{i}^{2}$, i.e., $J_{i} = \{j \in J: a_{ij} \ge b_{i}  \}$. \\
Based on Lemma 1 and Definition 2, we have also a necessary feasibility condition for Problem (2) as follows:\\

\textbf{Corollary 1.} If $S(\textbf{A}, \textbf{b}) \ne \emptyset$, then $J_{i} \ne \emptyset$, $\forall i \in I$. \\
\textbf{Proof.} Let $S(\textbf{A}, \textbf{b}) \ne \emptyset$. By contradiction, suppose that there exists some $i_{0} \in I$ such that $J_{i_{0}} = \emptyset$. So, according to Definition 2, for each $j \in J$ we have $a_{i_{0}j} < b_{i_{0}}$, and therefore $\min\{a_{i_{0}j}, x_{i_{0}}, x_{j} \} < b_{i_{0}}$ for each $j \in J$ and each $\textbf{x} \in [0,1]^{n}$. Consequently, from Lemma 1 we conclude that $S(\textbf{a}_{i},b_{i}) = \emptyset$ which contradicts $S(\textbf{A}, \textbf{b}) \ne \emptyset$. \\

In contrast to the feasible solutions set of Problem (1) (depicted in Figure 1), it will be proved that the feasible region of (3) (also, that of Problem (2)) interestingly comes in three different shapes depending on the value of the diagonal element $a_{ii}$. Figure 2 below schematically shows these differences in three cases $a_{ii} > b_{i}$, $a_{ii} = b_{i}$, and $a_{ii} < b_{i}$. \\

\begin{figure}[ht]
    \begin{center}
	    \includegraphics[height=8cm]{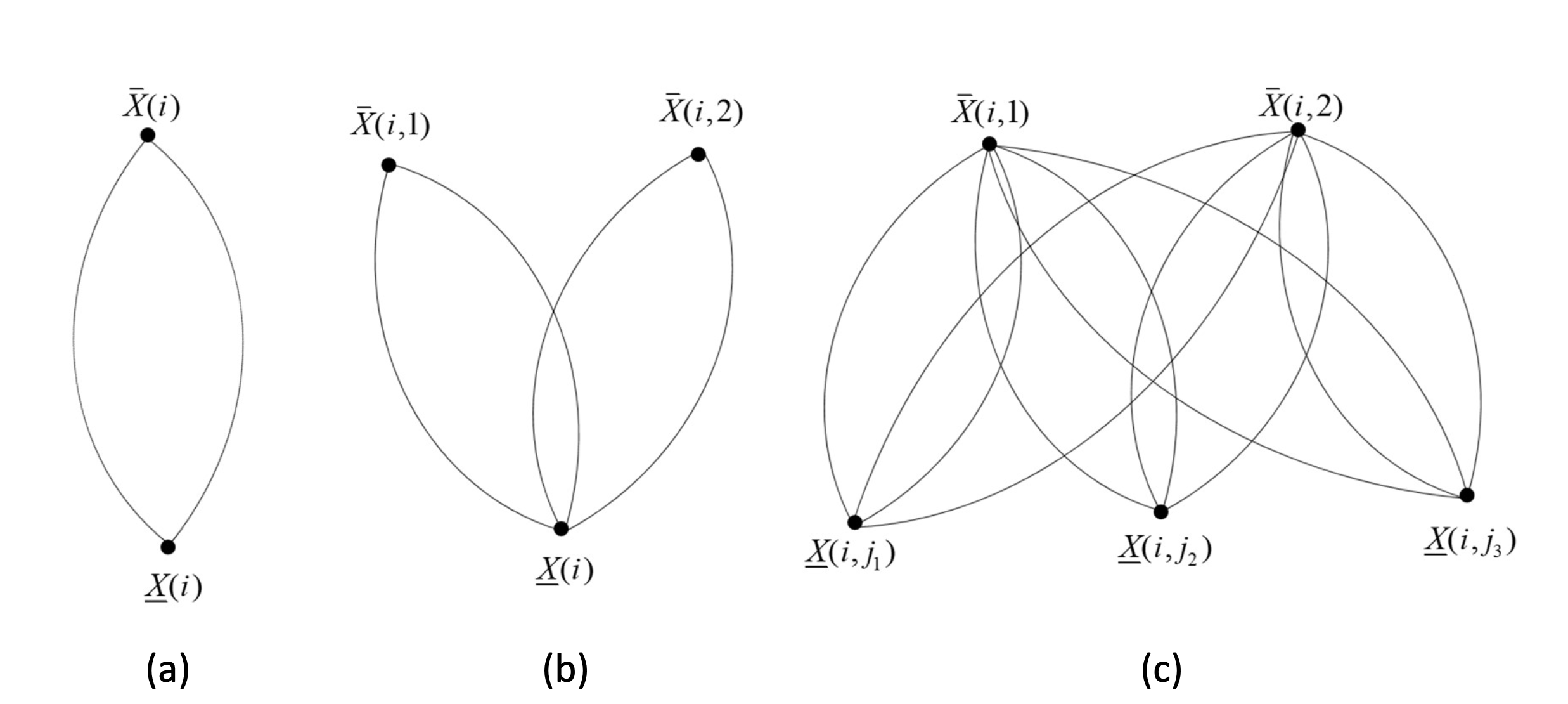}
	    \caption{The shape of the feasible solutions set for equation $\max_{j \in J}\{\min\{a_{ij}, x_{i}, x_{j} \}\} = b_{i}$, (a) if $a_{ii} > b_{i}$, (b) if $a_{ii} = b_{i}$, and (c) if $a_{ii} < b_{i}$.}
	\end{center}
\end{figure} 

As shown in Figure 2, when $a_{ii} > b_{i}$, set $S(\textbf{a}_{i},b_{i})$ has exactly one maximum solution, $\Bar{\textbf{X}}(i)$, and one minimum solution, $\underline{\textbf{X}}(i)$. In the case that  $a_{ii} = b_{i}$, $S(\textbf{a}_{i},b_{i})$ has yet a unique minimum, but exactly two maximal solutions $\Bar{\textbf{X}}(i,1)$ and $\Bar{\textbf{X}}(i,2)$. Finally, if $a_{ii} < b_{i}$, set $S(\textbf{a}_{i},b_{i})$ is determined by two maximal solutions and a finite number of minimal ones. It will be shown that the maximal solutions are always constructed in two specific ways, whether $a_{ii} = b_{i}$ or $a_{ii} < b_{i}$. For this reason, we interpret these two ways (possibilities) as two types for maximal solutions and this is why they are denoted by $\Bar{\textbf{X}}(i,1)$ and $\Bar{\textbf{X}}(i,2)$. So, solution $\Bar{\textbf{X}}(i,1)$ (solution $\Bar{\textbf{X}}(i,2)$) means the maximal solution type 1 (type 2) obtained by the first (second) way. As mentioned above, in the case that $a_{ii} \le b_{i}$, $S(\textbf{a}_{i},b_{i})$ includes exactly two maximal solutions $\Bar{\textbf{X}}(i,1)$ and $\Bar{\textbf{X}}(i,2)$. Moreover, it will be slightly later shown that if $a_{ii} > b_{i}$, the unique maximum solution $\Bar{\textbf{X}}(i)$ is also obtained by the same way used for constructing the maximal solutions type 1. In order to prove the correctness of the foregoing results, we firstly give some necessary definitions. \\

\textbf{Definition 3.} Let $I_{1} = \{i \in I : a_{ii} > b_{i} \}$, $I_{2} = \{i \in I : a_{ii} = b_{i} \}$, and $I_{3} = \{i \in I : a_{ii} < b_{i} \}$. \\

\textbf{Lemma 2.} Let $\textbf{x} \in S(\textbf{a}_{i},b_{i})$. \\
(a) $x_{i} \ge b_{i}$ \\
(b) If $i \in I_{1} \bigcup I_{2}$ and $x_{i} = b_{i}$, then $x_{j} \in [0,1], \forall j \in J - \{i \}$.\\
(c) If $i \in I_{3}$ and $x_{i} = b_{i}$, then there exists at least one $j_{0} \in J_{i}$ such that $x_{j_{0}} \ge b_{i}$.\\
(d) If $i \in I_{1}$, then $x_{i} = b_{i}$.\\
(e) If If $i \in I_{2} \bigcup I_{3}$ and $x_{i} > b_{i}$, then $x_{j} \le b_{i}$ for each $j \in J^{1}_{i}$.\\
(f) If $i \in I_{3}$, and $x_{i} > b_{i}$, then there exist some $j_{0} \in J_{i}$ such that either $j_{0} \in J^{1}_{i}$ and $x_{j_{0}} = b_{i}$ or $j_{0} \in J^{2}_{i}$ and $x_{j_{0}} \ge b_{i}$\\

\textbf{Proof.} By contradiction, suppose that $x_{i} < b_{i}$. So, $\min\{a_{ij}, x_{i}, x_{j}\} < b_{i}, \forall j \in J$. Hence, $\textbf{x}$ violates Part (II) of the necessary and sufficient feasibility conditions (4) and therefore $\textbf{x} \notin S(\textbf{a}_{i},b_{i})$ that is a contradiction. \textbf{(b)} By the assumptions, we have $\min\{a_{ii}, x_{i}, x_{i}\} = \min\{a_{ii}, b_{i}, b_{i}\} = b_{i}$ and $\min\{a_{ij}, x_{i}, x_{j}\} = \min\{a_{ij}, b_{i}, x_{j}\} \le b_{i}$, for each $j \in J - \{i\}$ and each $x_{j} \in [0,1]$. Hence, conditions (4) hold true for each value $x_{j} \in [0,1]$ where $J - \{i\}$. \textbf{(c)} The result is obtained by contradiction, i.e., by assuming that $x_{j} < b_{i}, \forall j \in J_{i}$. Hence, $\min\{a_{ij}, x_{i}, x_{j}\} < b_{i}, \forall j \in J_{i}$. On the other hand, we have $\min\{a_{ii}, x_{i}, x_{i}\} < b_{i}$ (because, $a_{ii} < b_{i}$ if $i \in I_{3}$) and 
and $\min \left\{a_{i j}, x_{i}, x_{j}\right\}<b_{i}, \forall j \notin J_{i}$ (because, $a_{i j}<b_{i}$ if $j \notin J_{i}$ ). Consequently, $\min \left\{a_{i j}, x_{i}, x_{j}\right\}<b_{i}, \forall j \in J$, that contradicts Part (II) of conditions (4). (d) According to Part (a), $x_{i} \geq b_{i}$. If $x_{i}>b_{i}$, then $\min \left\{a_{i i}, x_{i}, x_{i}\right\}>b_{i}$ which contradicts Part (I) of conditions (4). Therefore, we must have $x_{i}=b_{i}$. (e) By contradiction, suppose that $x_{j_{0}}>b_{i}$ for some $j_{0} \in J_{i}^{1}$. Thus, $\min \left\{a_{i_{0}}, x_{i}, x_{j_{0}}\right\}>b_{i}$ which violates Part (I) of conditions (4). (f) Since $i \in I_{3}$, then $a_{i i}<b_{i}$ that implies $\min \left\{a_{i i}, x_{i}, x_{i}\right\}<b_{i}$. By contradiction, suppose that $x_{j} \neq b_{i}, \forall j \in J_{i}^{1}$, and $x_{j}<b_{i}, \forall j \in J_{i}^{2}$. Therefore, $\min \left\{a_{i j}, x_{i}, x_{j}\right\}<b_{i}, \forall j \in J_{i}^{2}$. Now, if $x_{j}<b_{i}, \forall j \in J_{i}^{1}$, then we also have $\min \left\{a_{i j}, x_{i}, x_{j}\right\}<b_{i}, \forall j \in J_{i}^{1}$. So, $\min \left\{a_{i j}, x_{i}, x_{j}\right\}<b_{i}, \forall j \in J$, that contradicts Part (II) of conditions (4). Otherwise, if $x_{j_{0}}>b_{i}$ for some $j_{0} \in J_{i}^{1}$, then $\min \left\{a_{ij_{0}}, x_{i}, x_{j_{0}}\right\}>b_{i}$ which violates Part (I) of conditions (4).\\

\textbf{Definition 4.} For each $i \in I_{1}$, define two $n \times 1$ vectors \\ $\bar{X}(i)^{T}=\left(\bar{X}(i)_{1}, \bar{X}(i)_{2}, \ldots, \bar{X}(i)_{n}\right)$ and $\underline{\boldsymbol{X}}(i)^{T}=\left(\underline{X}(i)_{1}, \underline{X}(i)_{2}, \ldots, \underline{X}(i)_{n}\right)$ such that

$$
\bar{X}(i)_{j}=\left\{\begin{array}{ll}
b_{i} & j=i \\
1 & j \neq i
\end{array}, \quad \underline{X}(i)_{j}= \begin{cases}b_{i} & j=i \\
0 & j \neq i\end{cases}\right.
$$

for each $j \in J$.\\

\textbf{Definition 5.} For each $i \in I_{2}$, define three $n \times 1$ vectors
$$\begin{array}{l}
 \bar{X}(i, 1)^{T}=\left(\bar{X}(i, 1)_{1}, \bar{X}(i, 1)_{2}, \ldots, \bar{X}(i, 1)_{n}\right),\\
 \overline{\boldsymbol{X}}(i, 2)^{T}=\left(\bar{X}(i, 2)_{1}, \bar{X}(i, 2)_{2}, \ldots, \bar{X}(i, 2)_{n}\right), \\
 \underline{\boldsymbol{X}}(i)^{T}=\left(\underline{X}(i)_{1}, \underline{X}(i)_{2}, \ldots, \underline{X}(i)_{n}\right) 
\end{array}$$
such that

$$
\bar{X}(i, 1)_{j}=\left\{\begin{array}{cc}
b_{i} & j=i \\
1 & j \neq i
\end{array} \quad, \quad \bar{X}(i, 2)_{j}=\left\{\begin{array}{ll}
b_{i} & a_{i j}>b_{i} \\
1 & a_{i j} \leq b_{i}
\end{array} \quad, \quad \underline{X}(i)_{j}= \begin{cases}b_{i} & j=i \\
0 & j \neq i\end{cases}\right.\right.
$$
for each $j \in J$.\\

\textbf{Definition 6.} For each $i \in I_{3}$, two $n \times 1$ vectors $\bar{X}(i, 1)$ and $\bar{X}(i, 2)$ are defined as in Definition 5. Moreover, for each $i \in I_{3}$ and each $j \in J_{i}$, we define an $n \times 1$ vector $\underline{X}(i, j)^{T}=\left(\underline{X}(i, j)_{1}, \underline{X}(i, j)_{2}, \ldots, \underline{X}(i, j)_{n}\right)$ such that

$$
\underline{X}(i, j)_{k}= \begin{cases}b_{i} & k=i \text { or } k=j \\ 0 & \text { otherwise }\end{cases}
$$

for each $k \in J$.\\

Lemma 3 describes the shape of the feasible solutions set of $S\left(\boldsymbol{a}_{i}, b_{i}\right)$ for each value of the diagonal element $a_{i i}$, or equivalently for each $i \in I\left(I=I_{1} \cup I_{2} \cup I_{3}\right)$. In all the cases, notation $[\boldsymbol{Y}, \boldsymbol{Z}]$ (where $\boldsymbol{Y}=\left(Y_{j}\right)_{n \times 1}$ and $\boldsymbol{Z}=\left(Z_{j}\right)_{n \times 1}$ are two $n \times 1$ vectors) means the closed convex cell including all the $n \times 1$ vectors $\boldsymbol{x}=\left(x_{j}\right)_{n \times 1}$ such that $\boldsymbol{Y} \leq \boldsymbol{x} \leq \boldsymbol{Z}$, i.e., $Y_{j} \leq x_{j} \leq Z_{j}$ for each $j \in J$.\\

\textbf{Lemma 3.} \textbf{(a)} If $i \in I_{1}$, then $S\left(\boldsymbol{a}_{i}, b_{i}\right)=[\underline{X}(i), \bar{X}(i)]$.\\
\textbf{(b)} If $i \in I_{2}$, then $S\left(\boldsymbol{a}_{i}, b_{i}\right)=[\underline{\boldsymbol{X}}(i), \overline{\boldsymbol{X}}(i, 1)] \cup[\underline{\boldsymbol{X}}(i), \overline{\boldsymbol{X}}(i, 2)]$.\\
\textbf{(c)} If $i \in I_{3}$, then $S\left(\boldsymbol{a}_{i}, b_{i}\right)=\left\{\bigcup_{j \in J_{i}}[\underline{\boldsymbol{X}}(i, j), \overline{\boldsymbol{X}}(i, 1)]\right\} \bigcup\left\{\bigcup_{j \in J_{i}}[\underline{\boldsymbol{X}}(i, j), \overline{\boldsymbol{X}}(i, 2)]\right\}$.\\

\textbf{Proof.} (a) Since $i \in I_{1}$, then $a_{i i}>b_{i}$. Let $\boldsymbol{x} \in[\underline{X}(i), \overline{\boldsymbol{X}}(i)]$. From Definition 4, we have $\underline{X}(i)_{i}=x_{i}=\bar{X}(i)_{i}=b_{i}$. So, $\min \left\{a_{i i}, x_{i}, x_{i}\right\}=\min \left\{a_{i i}, b_{i}, b_{i}\right\}=b_{i}$ and $\min \left\{a_{i j}, x_{i}, x_{j}\right\}=\min \left\{a_{i j}, b_{i}, x_{j}\right\} \leq b_{i}$ for each $j \in J$ such that $j \neq i$. Hence, $\max _{j \in J}\left\{\min \left\{a_{i j}, x_{i}, x_{j}\right\}\right\}=b_{i}$ that implies $\boldsymbol{x} \in S\left(\boldsymbol{a}_{i}, b_{i}\right)$. Conversely, let $\boldsymbol{x} \in S\left(\boldsymbol{a}_{i}, b_{i}\right)$. So, from Lemma 2(d), we have $x_{i}=b_{i}$, and then the result follows from Lemma 2(b) and Definition 4. (b) By the assumption, $i \in I_{2}$ and therefore $a_{i i}=b_{i}$. Let $\boldsymbol{x} \in[\underline{\boldsymbol{X}}(i), \overline{\boldsymbol{X}}(i, 1)] \cup[\underline{\boldsymbol{X}}(i), \overline{\boldsymbol{X}}(i, 2)]$. If $\boldsymbol{x} \in[\underline{\boldsymbol{X}}(i), \overline{\boldsymbol{X}}(i, 1)]$, then the argument stated in Part (a) results in $\boldsymbol{x} \in S\left(\boldsymbol{a}_{i}, b_{i}\right)$. Otherwise, suppose $\boldsymbol{x} \in[\underline{X}(i), \overline{\boldsymbol{X}}(i, 2)]$ . Then, by Definition 5 we have $b_{i} \leq x_{i} \leq 1,0 \leq x_{j} \leq b_{i}$ when $j \in J_{i}^{1}$ and $0 \leq x_{j} \leq 1$ when $j \in J_{i}^{2}-\{i\}$ or $j \notin J_{i}$. Hence, $\min \left\{a_{i i}, x_{i}, x_{i}\right\}=b_{i}, \min \left\{a_{i j}, x_{i}, x_{j}\right\} \leq b_{i}$ if $j \in J_{i}-\{i\}$, and $\min \left\{a_{i j}, x_{i}, x_{j}\right\}<b_{i}$ if $j \notin J_{i}$ (because, $a_{i i}<b_{i}$ if $\left.j \notin J_{i}\right)$. Consequently, $\max _{j \in J}\left\{\min \left\{a_{i j}, x_{i}, x_{j}\right\}\right\}=b_{i}$ , i.e., $\boldsymbol{x} \in S\left(\boldsymbol{a}_{i}, b_{i}\right)$. The converse statement is resulted from Lemma 2 (Parts (a), (b) and (e)) and Definition 5. (c) Let $i \in I_{3}$ and $\boldsymbol{x} \in\left\{\bigcup_{j \in J_{i}}[\underline{\boldsymbol{X}}(i, j), \overline{\boldsymbol{X}}(i, 1)]\right\} \cup\left\{\bigcup_{j \in J_{i}}[\underline{\boldsymbol{X}}(i, j), \overline{\boldsymbol{X}}(i, 2)]\right\}$. Therefore, $\quad a_{i i}<b_{i} . \quad$ If $\boldsymbol{x} \in \bigcup_{j \in J_{i}}[\underline{\boldsymbol{X}}(i, j), \overline{\boldsymbol{X}}(i, 1)]$, then there exists at least one $j^{\prime} \in J_{i}$ such that $\boldsymbol{x} \in\left[\underline{\boldsymbol{X}}\left(i, j^{\prime}\right), \overline{\boldsymbol{X}}(i, 1)\right]$ . Thus, by Definition 6, $x_{i}=b_{i}$ and $b_{i} \leq x_{j^{\prime}} \leq 1$. Therefore, we have $\min \left\{a_{i i}, x_{i}, x_{i}\right\}<b_{i}$, $\min \left\{a_{i^{\prime}}, x_{i}, x_{j^{\prime}}\right\}=b_{i}, \quad \min \left\{a_{i j}, x_{i}, x_{j}\right\} \leq b_{i}, \quad \forall j \in J_{i}-\left\{j^{\prime}\right\}$, and $\min \left\{a_{i j}, x_{i}, x_{j}\right\}<b_{i}$ when $j \notin J_{i}$. Hence, $\max _{j \in J}\left\{\min \left\{a_{i j}, x_{i}, x_{j}\right\}\right\}=b_{i}$ that means $\boldsymbol{x} \in S\left(\boldsymbol{a}_{i}, b_{i}\right)$. On the other hand, if $\boldsymbol{x} \in \bigcup_{j \in J_{i}}[\underline{\boldsymbol{X}}(i, j), \overline{\boldsymbol{X}}(i, 2)]$, then $\boldsymbol{x} \in\left[\underline{\boldsymbol{X}}\left(i, j^{\prime}\right), \overline{\boldsymbol{X}}(i, 2)\right]$ for at least one $j^{\prime} \in J_{i}$. In this case, Definition 6 implies $b_{i} \leq x_{i} \leq 1$; also, we have $x_{j^{\prime}}=b_{i}$ if $j^{\prime} \in J_{i}^{1}$ and $b_{i} \leq x_{j^{\prime}} \leq 1$ if $j^{\prime} \in J_{i}^{2}$; additionally, $0 \leq x_{j} \leq b_{i}$ if $j \in J_{i}^{1}-\left\{j^{\prime}\right\}$ and $0 \leq x_{j} \leq 1$ if $j \in J_{i}^{2}-\left\{j^{\prime}\right\}$. Therefore, $\min \left\{a_{i i}, x_{i}, x_{i}\right\}<b_{i}, \min \left\{a_{i j^{\prime}}, x_{i}, x_{j^{\prime}}\right\}=b_{i}, \min \left\{a_{i j}, x_{i}, x_{j}\right\} \leq b_{i}, \forall j \in J_{i}-\left\{j^{\prime}\right\}$, and $\min \left\{a_{i j}, x_{i}, x_{j}\right\}<b_{i}$ if $j \notin J_{i}$. So, we obtain the same result, i.e., $\max _{j \in J}\left\{\min \left\{a_{i j}, x_{i}, x_{j}\right\}\right\}=b_{i}$. Conversely, let $\boldsymbol{x} \in S\left(\boldsymbol{a}_{i}, b_{i}\right)$. From Lemma 2(a), $x_{i} \geq b_{i}$. If $x_{i}=b_{i}$, then by Lemma 2(c) and Definition 6, $\boldsymbol{x} \in\left[\underline{X}\left(i, j_{0}\right), \overline{\boldsymbol{X}}(i, 1)\right]$ for some $j_{0} \in J_{i}$ such that $x_{j_{0}} \geq b_{i}$. Otherwise, if $x_{i}>b_{i}$, then Lemma 2(e) and Definition 6 imply $\boldsymbol{x} \leq \overline{\boldsymbol{X}}(i, 2)$. Moreover, from Lemma 2(f) and Definition 6 , we have $\underline{X}\left(i, j_{0}\right) \leq \boldsymbol{x}$ for some $j_{0} \in J_{i}$ such that either $j_{0} \in J_{i}^{1}$ and $x_{j_{0}}=b_{i}$ or $j_{0} \in J_{i}^{2}$ and $x_{j_{0}} \geq b_{i}$. As a result, we have $\boldsymbol{x} \in\left[\underline{\boldsymbol{X}}\left(i, j_{0}\right), \overline{\boldsymbol{X}}(i, 2)\right]$ where $j_{0} \in J_{i}$.\\

\textbf{Corollary 2.} According to Lemma 3, solutions $\bar{X}(i)$ and $\underline{X}(i)$ are the maximum and minimum solutions of $S\left(\boldsymbol{a}_{i}, b_{i}\right), \forall i \in I_{1}$, respectively. Also, there are two maximal solutions $\overline{\boldsymbol{X}}(i, 1)$ and $\overline{\boldsymbol{X}}(i, 2)$, and a unique minimum solution $\underline{\boldsymbol{X}}(i)$ for $S\left(\boldsymbol{a}_{i}, b_{i}\right), \forall i \in I_{2}$. Moreover, for each $i \in I_{3}, \quad S\left(\boldsymbol{a}_{i}, b_{i}\right)$ has two maximal solutions $\bar{X}(i, 1)$ and $\bar{X}(i, 2)$, and a finite number of minimal solutions $\underline{X}(i, j), \forall j \in J_{i}$. Additionally, from Definitions 46 , it is clear that the maximum solution $\bar{X}(i)$ is defined in the same way as the maximal solutions $\overline{\boldsymbol{X}}(i, 1)$ are defined. Lemma 3 determines the feasible region of the single equation $\boldsymbol{a}_{i} \otimes \boldsymbol{x}=b_{i}$ in all the cases, where the value of $a_{i i}$ may be greater than $b_{i}$, less than $b_{i}$ or equal to $b_{i}$. Based on the differences between the feasible regions of these single equations depending on the values of $a_{i i}$, it seems reasonable to categorize the equations of Problem (2) into three groups of constraints such that $k^{\prime}$th group includes each equation $\boldsymbol{a}_{i} \otimes \boldsymbol{x}=b_{i}$ (i.e., $\max _{j \in J}\left\{\min \left\{a_{i j}, x_{i}, x_{j}\right\}\right\}=b_{i}$ ) with $i \in I_{k}$, for $k=1,2,3$. The following definition formally expresses the feasible regions of these three groups of the constraints. In what follows, we separately investigate these feasible regions and determine their shapes as well as their extreme solutions (maximum, minimum, maximal and minimal solutions). Obviously, the intersection of these feasible regions will finally give the feasible solutions set of Problem (2).\\

\textbf{Definition 7.} Let $S_{k}(\boldsymbol{A}, \boldsymbol{b})=\bigcap_{i \in I_{k}} S\left(\boldsymbol{a}_{i}, b_{i}\right), k=1,2,3$.\\

Based on the notation used in Definition 7, $S(\boldsymbol{A}, \boldsymbol{b})=S_{1}(\boldsymbol{A}, \boldsymbol{b}) \cap S_{2}(\boldsymbol{A}, \boldsymbol{b}) \cap S_{3}(\boldsymbol{A}, \boldsymbol{b})$. The following definition gives a useful tool enabling us to determine the minimal and maximal solutions of sets $S_{2}(\boldsymbol{A}, \boldsymbol{b})$ and $S_{3}(\boldsymbol{A}, \boldsymbol{b})$.\\

\textbf{Definition 8.} Let $e^{\prime}: I_{2} \rightarrow\{1,2\}$ be a function from $I_{2}$ into $\{1,2\}$ and $E^{\prime}$ denote the set of all the functions $e^{\prime}$ on $I_{2}$. Similarly, the notation $E^{\prime \prime}$ is used to denote the set of all functions $e^{\prime \prime}: I_{3} \rightarrow\{1,2\}$. Moreover, let $\underline{e}: I_{3} \rightarrow \bigcup_{i \in I_{3}} J_{i}$ be a function from $I_{3}$ into $\bigcup_{i \in I_{3}} J_{i}$ such that $\underline{e}(i) \in J_{i}, \forall i \in I_{3}$, and let $\underline{E}$ denote the set of all the functions $\underline{e}$. For the sake of simplicity, each $\underline{e} \in \underline{E}$ can be also presented as the vector $\underline{e}=\left[j_{1}, j_{2}, \ldots, j_{s}\right]$ in which $j_{k}=\underline{e}(k), k=1,2, \ldots, s$.\\

\textbf{Remark 1.} According to Definition 8, for each $e^{\prime} \in E^{\prime}$ and each $i \in I_{2}, e^{\prime}(i)$ can be interpreted as a variable that takes values in $\operatorname{dom}\left(e^{\prime}(i)\right)=\{1,2\}$. By the same interpretation, $\operatorname{dom}\left(e^{\prime \prime}(i)\right)=\{1,2\}$ for each $e^{\prime \prime} \in E^{\prime \prime}$ and each $i \in I_{3}$. Similarly, for each $\underline{e} \in \underline{E}$ and each $i \in I_{3}$, variable $\underline{e}(i)$ takes values in $\operatorname{dom}(\underline{e}(i))=J_{i}$. Furthermore, $\left|E^{\prime}\right|=2^{\left|I_{2}\right|}$ , $\left|E^{\prime \prime}\right|=2^{\left|I_{3}\right|}$ and $|\underline{E}|=\prod_{i \in I_{3}}\left|J_{i}\right|$, where $|\cdot|$ denotes the cardinality of the sets. \\

\textbf{Definition 9.} Define $\underline{\boldsymbol{X}}_{1}=\max _{i \in I_{1}}\{\underline{\boldsymbol{X}}(i)\}$ and $\overline{\boldsymbol{X}}_{1}=\min _{i \in I_{1}}\{\overline{\boldsymbol{X}}(i)\}$. Also, let $\underline{\boldsymbol{X}}_{2}=\max _{i \in I_{2}}\{\underline{\boldsymbol{X}}(i)\}$ and $\overline{\boldsymbol{X}}_{2}\left(e^{\prime}\right)=\min _{i \in I_{2}}\left\{\overline{\boldsymbol{X}}\left(i, e^{\prime}(i)\right)\right\}, \forall e^{\prime} \in E^{\prime}$. Additionally, for each $\underline{e} \in \underline{E}$ and each $e^{\prime \prime} \in E^{\prime \prime}$, we define $\underline{\boldsymbol{X}}_{3}(\underline{e})=\max _{i \in I_{3}}\{\underline{\boldsymbol{X}}(i, \underline{e}(i))\}$ and $\overline{\boldsymbol{X}}_{3}\left(e^{\prime \prime}\right)=\min _{i \in I_{3}}\left\{\overline{\boldsymbol{X}}\left(i, e^{\prime \prime}(i)\right)\right\}$, respectively.

The following lemma shows that the vectors defined in Definition 9 are indeed the extreme solutions determining the feasible regions of sets $S_{1}(\boldsymbol{A}, \boldsymbol{b}), S_{2}(\boldsymbol{A}, \boldsymbol{b})$ and $S_{3}(\boldsymbol{A}, \boldsymbol{b})$.\\

\textbf{Lemma 4.} \textbf{(a)} $S_{1}(\boldsymbol{A}, \boldsymbol{b})=\left[\underline{X}_{1}, \overline{\boldsymbol{X}}_{1}\right]$.

\textbf{(b)} $S_{2}(\boldsymbol{A}, \boldsymbol{b})=\bigcup_{e^{\prime} \in E^{\prime}}\left[\underline{\boldsymbol{X}}_{2}, \overline{\boldsymbol{X}}_{2}\left(e^{\prime}\right)\right]$.

\textbf{(c)} $S_{3}(\boldsymbol{A}, \boldsymbol{b})=\bigcup_{\underline{e} \in \underline{E}} \bigcup_{ e^{\prime \prime}  \in E^{\prime \prime}}\left[\underline{\boldsymbol{X}}_{3}(\underline{e}), \overline{\boldsymbol{X}}_{3}\left(e^{\prime \prime}\right)\right]$.

\textbf{Proof.} \textbf{(a)} From Lemma 3(a) and Definitions 7 and 9, we have

$$
S_{1}(\boldsymbol{A}, \boldsymbol{b})=\bigcap_{i \in I_{1}}[\underline{\boldsymbol{X}}(i), \overline{\boldsymbol{X}}(i)]=\left[\max _{i \in I_{1}}\{\underline{\boldsymbol{X}}(i)\}, \min _{i \in I_{1}}\{\overline{\boldsymbol{X}}(i)\}\right]=\left[\underline{\boldsymbol{X}}_{1}, \overline{\boldsymbol{X}}_{1}\right]
$$

\textbf{(b)} The result is obtained from the following equalities, Lemma 3(b) and Definitions 7-9:

$$
\begin{aligned}
S_{2}(\boldsymbol{A}, \boldsymbol{b}) & =\bigcap_{i \in I_{2}}\{[\underline{\boldsymbol{X}}(i), \overline{\boldsymbol{X}}(i, 1)] \cup[\underline{\boldsymbol{X}}(i), \overline{\boldsymbol{X}}(i, 2)]\}=\bigcup_{e^{\prime} \in E^{\prime}}\left\{\bigcap_{i \in I_{2}}\left[\underline{\boldsymbol{X}}(i), \overline{\boldsymbol{X}}\left(i, e^{\prime}(i)\right)\right]\right\} \\
& =\bigcup_{e^{\prime} \in E^{\prime}}\left[\max _{i \in I_{2}}\{\underline{\boldsymbol{X}}(i)\}, \min _{i \in I_{2}}\left\{\overline{\boldsymbol{X}}\left(i, e^{\prime}(i)\right)\right\}\right]=\bigcup_{e^{\prime} \in E^{\prime}}\left[\underline{\boldsymbol{X}}_{2}, \overline{\boldsymbol{X}}_{2}\left(e^{\prime}\right)\right]
\end{aligned}
$$\\

\textbf{(c)} Similarly, by Lemma 3(c) and Definitions 7-9, we have

$$
\begin{aligned}[t]
& S_{3}(\boldsymbol{A}, \boldsymbol{b})=\bigcap_{i \in I_{3}}\left\{\left\{\bigcup_{j \in J_{i}}[\underline{\boldsymbol{X}}(i, j), \overline{\boldsymbol{X}}(i, 1)]\right\} \bigcup\left\{\bigcup_{j \in J_{i}}[\underline{\boldsymbol{X}}(i, j), \overline{\boldsymbol{X}}(i, 2)]\right\}\right\} \\
& =\bigcup_{e^{\prime \prime} \in E^{\prime \prime}}\left\{\bigcap_{i \in I_{3}}\left\{\bigcup_{j \in J_{i}}\left[\underline{\boldsymbol{X}}(i, j), \overline{\boldsymbol{X}}\left(i, e^{\prime \prime}(i)\right)\right]\right\}\right\} \\
& = \bigcup_{e^{\prime \prime} \in E^{\prime \prime}}\left\{\bigcup_{\underline e \in \underline E}\left\{\bigcap_{i \in I_{3}}\left[\underline{\boldsymbol{X}}(i, \underline e(i)), \overline{\boldsymbol{X}}\left(i, e^{\prime \prime}(i)\right)\right]\right\}\right\} \\
& = \bigcup_{e^{\prime \prime} \in E^{\prime \prime}}\bigcup_{\underline e \in \underline E}\{\bigcap_{i \in I_{3}} \left[\underline{\boldsymbol{X}}(i, \underline e(i)), \overline{\boldsymbol{X}}\left(i, e^{\prime \prime}(i)\right)\right]\} \\
& = \bigcup_{e^{\prime \prime} \in E^{\prime \prime}} \bigcup_{\underline e \in \underline E}[\max_{i \in I_{3}}\{\underline{\boldsymbol{X}}(i, \underline e(i))\}, \min_{i \in I_{3}}\{ \overline{\boldsymbol{X}}(i, e^{\prime \prime}(i))\}] ~ = ~ \bigcup_{e^{\prime \prime} \in E^{\prime \prime}} \bigcup_{\underline e \in \underline E} [ \underline{\boldsymbol{X}}_{3}(\underline e), \overline{\boldsymbol{X}}_{3}(e^{\prime \prime}) ]
\end{aligned}
$$

Thus Part (c) holds, and the proof is complete.

Based on Lemma 4(b), $S_{2}(\boldsymbol{A}, \boldsymbol{b})$ is formed as the union of a finite number of nonempty cells $\left[\underline{\boldsymbol{X}}_{2}, \overline{\boldsymbol{X}}_{2}\left(e^{\prime}\right)\right]\left(e^{\prime} \in E^{\prime}\right)$. Clearly, $\left[\underline{\boldsymbol{X}}_{2}, \overline{\boldsymbol{X}}_{2}\left(e^{\prime}\right)\right] \neq \varnothing$ iff $\underline{\boldsymbol{X}}_{2} \leq \overline{\boldsymbol{X}}_{2}\left(e^{\prime}\right)$. So, $\underline{\boldsymbol{X}}_{2}$ is the unique minimum solution of $S_{2}(\boldsymbol{A}, \boldsymbol{b})$, and also each $\overline{\boldsymbol{X}}_{2}\left(e^{\prime}\right) \quad\left(e^{\prime} \in E^{\prime}\right)$ such that $\underline{\boldsymbol{X}}_{2} \leq \overline{\boldsymbol{X}}_{2}\left(e^{\prime}\right)$ is a maximal solution of $S_{2}(\boldsymbol{A}, \boldsymbol{b})$. Similarly, from Lemma 4(c), if $\left[\underline{\boldsymbol{X}}_{3}(\underline{e}), \overline{\boldsymbol{X}}_{3}\left(e^{\prime \prime}\right)\right] \neq \varnothing$ for any $\underline{e} \in \underline{E}$ and $e^{\prime \prime} \in E^{\prime \prime}$, then $\underline{\boldsymbol{X}}_{3}(\underline{e})$ and $\overline{\boldsymbol{X}}_{3}\left(e^{\prime \prime}\right)$ are minimal and maximal solutions of $S_{3}(\boldsymbol{A}, \boldsymbol{b})$, respectively. From Parts (a) - (c) of Lemma 4, the following necessary and sufficient conditions can be directly resulted for the feasibility of sets $S_{k}(\boldsymbol{A}, \boldsymbol{b}), k=1,2,3$: \\

\textbf{Corollary 3.} \textbf{(a)} $S_{1}(\boldsymbol{A}, \boldsymbol{b}) \neq \varnothing$ iff $\underline{\boldsymbol{X}}_{1} \leq \overline{\boldsymbol{X}}_{1}$.\\
\textbf{(b)} $S_{2}(\boldsymbol{A}, \boldsymbol{b}) \neq \varnothing$ iff there exists some $e^{\prime} \in E^{\prime}$ such that $\underline{\boldsymbol{X}}_{2} \leq \overline{\boldsymbol{X}}_{2}\left(e^{\prime}\right)$.\\
\textbf{(c)} $S_{3}(\boldsymbol{A}, \boldsymbol{b}) \neq \varnothing$ iff there exist some $\underline{e} \in \underline{E}$ and $e^{\prime \prime} \in E^{\prime \prime}$ such that $\underline{\boldsymbol{X}}_{3}(\underline{e}) \leq \overline{\boldsymbol{X}}_{3}\left(e^{\prime \prime}\right)$. \\

The following theorem provides a necessary condition for the feasibility of Problem (2) that is used in the next section, where some rules are introduced for reducing the problem. \\ 

\textbf{Theorem 1.} Suppose that $S(\boldsymbol{A}, \boldsymbol{b}) \neq \varnothing$. Then, $\max \left\{\underline{\boldsymbol{X}}_{1}, \underline{\boldsymbol{X}}_{2}\right\} \leq \overline{\boldsymbol{X}}_{1}$. \\
\textbf{Proof.} Let $\quad \boldsymbol{x} \in S(\boldsymbol{A}, \boldsymbol{b})$. Since $\quad S(\boldsymbol{A}, \boldsymbol{b})=S_{1}(\boldsymbol{A}, \boldsymbol{b}) \cap S_{2}(\boldsymbol{A}, \boldsymbol{b}) \cap S_{3}(\boldsymbol{A}, \boldsymbol{b})$, then $\boldsymbol{x} \in S_{1}(\boldsymbol{A}, \boldsymbol{b}) \cap S_{2}(\boldsymbol{A}, \boldsymbol{b})$. So, Lemma 4 implies $\boldsymbol{x} \in\left[\underline{\boldsymbol{X}}_{1}, \overline{\boldsymbol{X}}_{1}\right] \cap\left\{\bigcup_{e^{\prime} \in E^{\prime}}\left[\underline{\boldsymbol{X}}_{2}, \overline{\boldsymbol{X}}_{2}\left(e^{\prime}\right)\right]\right\}$, i.e., $\boldsymbol{x} \in \bigcup_{e^{\prime} \in E^{\prime}}\left\{\left[\underline{\boldsymbol{X}}_{1}, \overline{\boldsymbol{X}}_{1}\right] \cap\left[\underline{\boldsymbol{X}}_{2}, \overline{\boldsymbol{X}}_{2}\left(e^{\prime}\right)\right]\right\}$. Therefore, $\boldsymbol{x} \in\left[\underline{\boldsymbol{X}}_{1}, \overline{\boldsymbol{X}}_{1}\right] \cap\left[\underline{\boldsymbol{X}}_{2}, \overline{\boldsymbol{X}}_{2}\left(e^{\prime}\right)\right]$ for at least one $e^{\prime} \in E^{\prime}$. Subsequently, we can conclude that $\left[\underline{\boldsymbol{X}}_{1}, \overline{\boldsymbol{X}}_{1}\right] \neq \varnothing$ and $\left[\underline{\boldsymbol{X}}_{1}, \overline{\boldsymbol{X}}_{1}\right] \cap\left[\underline{\boldsymbol{X}}_{2}, \overline{\boldsymbol{X}}_{2}\left(e^{\prime}\right)\right] \neq \varnothing$. Thus, we must have $\underline{\boldsymbol{X}}_{1} \leq \overline{\boldsymbol{X}}_{1}$ and $\underline{\boldsymbol{X}}_{2} \leq \overline{\boldsymbol{X}}_{1}$, and then the result follows.

Theorem 2 determines the extreme solutions of Problem (2) that are also the bounds of the feasible region.\\

\textbf{Theorem 2.} Suppose that $S(\boldsymbol{A}, \boldsymbol{b}) \neq \varnothing$. Then, 
$$\begin{array}{l}
S(\boldsymbol{A}, \boldsymbol{b}) = \bigcup_{\underline e \in \underline E} \bigcup_{e^{\prime \prime} \in E^{\prime \prime}} \bigcup_{e^{\prime} \in E^{\prime}} [\max\{ \underline{\boldsymbol{X}}_{1}, \underline{\boldsymbol{X}}_{2}, \underline{\boldsymbol{X}}_{3}(e) \}, \min\{\overline{\boldsymbol{X}}_{1}, \overline{\boldsymbol{X}}_{2}(e^{\prime}), \overline{\boldsymbol{X}}_{3}(e^{\prime \prime}) \}].
\end{array}$$

\textbf{Proof.} Since $S(\boldsymbol{A}, \boldsymbol{b})=S_{1}(\boldsymbol{A}, \boldsymbol{b}) \cap S_{2}(\boldsymbol{A}, \boldsymbol{b}) \cap S_{3}(\boldsymbol{A}, \boldsymbol{b})$, Lemma 4 implies that

$$\begin{array}{l}
S(\boldsymbol{A}, \boldsymbol{b})=\left[\underline{\boldsymbol{X}}_{1}, \overline{\boldsymbol{X}}_{1}\right] \cap\left\{\bigcup_{e^{\prime} \in E^{\prime}}\left[\underline{\boldsymbol{X}}_{2}, \overline{\boldsymbol{X}}_{2}\left(e^{\prime}\right)\right]\right\} \cap\left\{\bigcup_{\underline{e} \in \underline{E}} \bigcup_{e^{\prime \prime } \in E^{\prime \prime } }\left[\underline{\boldsymbol{X}}_{3}(\underline{e}), \overline{\boldsymbol{X}}_{3}\left(e^{\prime \prime}\right)\right]\right\}.
\end{array}$$

Therefore, we have

$$
\begin{aligned}
S(A, b) & =\bigcup_{\underline{e} \in \underline{E}} \bigcup_{e^{\prime \prime} \in E^{\prime \prime}} \bigcup_{e^{\prime} \in E^{\prime}}\left\{\left[\underline{\boldsymbol{X}}_{1}, \overline{\boldsymbol{X}}_{1}\right] \cap\left[\underline{\boldsymbol{X}}_{2}, \overline{\boldsymbol{X}}_{2}\left(e^{\prime}\right)\right] \cap\left[\underline{\boldsymbol{X}}_{3}(\underline{e}), \overline{\boldsymbol{X}}_{3}\left(e^{\prime \prime}\right)\right]\right\} \\
& =\bigcup_{\underline{e} \in \underline{E}} \bigcup_{e^{\prime \prime} \in E^{\prime \prime}} \bigcup_{e^{\prime} \in E^{\prime}}\left[\max \left\{\underline{\boldsymbol{X}}_{1}, \underline{\boldsymbol{X}}_{2}, \underline{\boldsymbol{X}}_{3}(\underline{e})\right\}, \min \left\{\overline{\boldsymbol{X}}_{1}, \overline{\boldsymbol{X}}_{2}\left(e^{\prime}\right), \overline{\boldsymbol{X}}_{3}\left(e^{\prime \prime}\right)\right\}\right]
\end{aligned}
$$
that completes the proof.\\

Based on Theorem 2, if $\left[\max \left\{\underline{\boldsymbol{X}}_{1}, \underline{\boldsymbol{X}}_{2}, \underline{\boldsymbol{X}}_{3}(\underline{e})\right\}, \min \left\{\overline{\boldsymbol{X}}_{1}, \overline{\boldsymbol{X}}_{2}\left(e^{\prime}\right), \overline{\boldsymbol{X}}_{3}\left(e^{\prime \prime}\right)\right\}\right] \neq \varnothing$ for any $\underline{e} \in \underline{E}$ , $e^{\prime} \in E^{\prime}$ and $e^{\prime \prime} \in E^{\prime \prime}$, then
\begin{align*}
\max \left\{\underline{\boldsymbol{X}}_{1}, \underline{\boldsymbol{X}}_{2}, \underline{\boldsymbol{X}}_{3}(\underline{e})\right\}
\end{align*}
and
\begin{align*}
\min \left\{\overline{\boldsymbol{X}}_{1}, \overline{\boldsymbol{X}}_{2}\left(e^{\prime}\right), \overline{\boldsymbol{X}}_{3}\left(e^{\prime \prime}\right)\right\}
\end{align*}
are minimal and maximal solutions of $S(\boldsymbol{A}, \boldsymbol{b})$, respectively. Theorem 2 also gives the following necessary and sufficient condition for the feasibility of Problem (2): \\

\textbf{Theorem 3.} $S(\boldsymbol{A}, \boldsymbol{b}) \neq \varnothing$ if and only if there exist $\underline{e} \in \underline{E}, e^{\prime} \in E^{\prime}$ and $e^{\prime \prime} \in E^{\prime \prime}$ such that $\max \left\{\underline{\boldsymbol{X}}_{1}, \underline{\boldsymbol{X}}_{2}, \underline{\boldsymbol{X}}_{3}(e)\right\} \leq \min \left\{\overline{\boldsymbol{X}}_{1}, \overline{\boldsymbol{X}}_{2}\left(e^{\prime}\right), \overline{\boldsymbol{X}}_{3}\left(e^{\prime \prime}\right)\right\}$.\\

\textbf{Proof.} The proof is directly resulted from Theorem 2.

\section{Simplification rules}

Based on Theorem 2, the feasible solutions set of Problem (2) is completely determined by a finite number of non-empty closed convex cells. Furthermore, by Theorem 3, $S(\boldsymbol{A}, \boldsymbol{b}) \neq \varnothing$ if and only if 
\begin{align*}
\left[\max \left\{\underline{\boldsymbol{X}}_{1}, \underline{\boldsymbol{X}}_{2}, \underline{\boldsymbol{X}}_{3}(e)\right\}, \min \left\{\overline{\boldsymbol{X}}_{1}, \overline{\boldsymbol{X}}_{2}\left(e^{\prime}\right), \overline{\boldsymbol{X}}_{3}\left(e^{\prime \prime}\right)\right\}\right] \neq \varnothing
\end{align*}
for some $\underline{e} \in \underline{E}, e^{\prime} \in E^{\prime}$ and $e^{\prime \prime} \in E^{\prime \prime}$. For this reason, the resolution of $S(\boldsymbol{A}, \boldsymbol{b})$ can be accelerated by focusing only on those functions $\underline{e} \in \underline{E}, e^{\prime} \in E^{\prime}$ and $e^{\prime \prime} \in E^{\prime \prime}$ such that

\begin{equation}
    \begin{array}{l}
        \max \left\{\underline{\boldsymbol{X}}_{1}, \underline{\boldsymbol{X}}_{2}, \underline{\boldsymbol{X}}_{3}(e)\right\} \leq \min \left\{\overline{\boldsymbol{X}}_{1}, \overline{\boldsymbol{X}}_{2}\left(e^{\prime}\right), \overline{\boldsymbol{X}}_{3}\left(e^{\prime \prime}\right)\right\}
    \end{array}
\end{equation}

For the sake of expository reference, such functions are formally defined in the following definition:\\

\textbf{Definition 10.} Suppose that $S(\boldsymbol{A}, \boldsymbol{b}) \neq \varnothing$ (and therefore, $J_{i} \neq \varnothing, \forall i \in I$, from Corollary 1), $I_{2} \neq \varnothing$ and $I_{3} \neq \varnothing$. Triple $\left(\underline{e}, e^{\prime}, e^{\prime \prime}\right)$ of functions $\underline{e} \in \underline{E}, e^{\prime} \in E^{\prime}$ and $e^{\prime \prime} \in E^{\prime \prime}$ is called admissible, if $\underline{e}, e^{\prime}$ and $e^{\prime \prime}$ satisfy Relation (5). Also, the set of all the admissible triples is denoted by $T$, i.e.,
\begin{align*}
T = \Biggl\{&\left(\underline{e}, e^{\prime}, e^{\prime \prime}\right): \underline{e} \in \underline{E}, e^{\prime} \in E^{\prime}, e^{\prime \prime} \in E^{\prime \prime}, \\
&\max \Bigl\{\underline{\boldsymbol{X}}_{1}, \underline{\boldsymbol{X}}_{2}, \underline{\boldsymbol{X}}_{3}(\underline{e})\Bigr\} \\
&\qquad \leq \min \Bigl\{\overline{\boldsymbol{X}}_{1}, \overline{\boldsymbol{X}}_{2}\left(e^{\prime}\right), \overline{\boldsymbol{X}}_{3}\left(e^{\prime \prime}\right)\Bigr\} \Biggr\} .
\end{align*}
In contrast, an inadmissible triple is a triple $\left(\underline{e}, e^{\prime}, e^{\prime \prime}\right)$ that violates Relation (5) and therefore leads to an empty cell $\left[\max \left\{\underline{\boldsymbol{X}}_{1}, \underline{\boldsymbol{X}}_{2}, \underline{\boldsymbol{X}}_{3}(\underline{e})\right\}, \min \left\{\overline{\boldsymbol{X}}_{1}, \overline{\boldsymbol{X}}_{2}\left(e^{\prime}\right), \overline{\boldsymbol{X}}_{3}\left(e^{\prime \prime}\right)\right\}\right]$. Consequently, the inadmissible triples can be considered as redundant or irrelevant selections of functions $\underline{e}, e^{\prime}$ and $e^{\prime \prime}$ in the sense that if they are disregarded, the feasible region $S(\boldsymbol{A}, \boldsymbol{b})$ is not affected. Due to this fact, some simplification rules are presented that restrict the selection possibilities of triples $\left(e, e^{\prime}, e^{\prime \prime}\right)$ by removing the inadmissible triples as far as possible. These rules are initially applied to Problem (2) to reduce the problem before starting to solve it, i.e., before selecting any triple $\left(\underline{e}, e^{\prime}, e^{\prime \prime}\right)$ to construct solutions $\underline{\boldsymbol{X}}_{3}(\underline{e}), \overline{\boldsymbol{X}}_{2}\left(e^{\prime}\right)$ and $\overline{\boldsymbol{X}}_{3}\left(e^{\prime \prime}\right)$. \\

\textbf{Definition 11.} Suppose that $I_{2}=\left\{i_{1}^{\prime}, \ldots, i_{p}^{\prime}\right\}$ and $I_{3}=\left\{i_{1}^{\prime \prime}, \ldots, i_{q}^{\prime \prime}\right\}$. Let $\overline{\boldsymbol{M}}^{1}=\left(\bar{m}_{i j}^{1}\right)_{p \times n}$ and $\overline{\boldsymbol{M}}^{2}=\left(\bar{m}_{i j}^{2}\right)_{p \times n}$ be two $p \times n$ matrices, and $\overline{\boldsymbol{N}}^{1}=\left(\bar{n}_{i j}^{1}\right)_{q \times n}$ and $\overline{\boldsymbol{N}}^{2}=\left(\bar{n}_{i j}^{2}\right)_{q \times n}$ be two $q \times n$ matrices whose $k$ 'th rows (denoted by $\overline{\boldsymbol{M}}_{k}^{1}, \overline{\boldsymbol{M}}_{k}^{2}, \bar{N}_{k}^{1}$ and $\bar{N}_{k}^{2}$, respectively) are defined as follows:

$$
\begin{aligned}
\overline{\boldsymbol{M}}_{\boldsymbol{k}}^{1}=\overline{\boldsymbol{X}}\left(i_{k}^{\prime}, 1\right), \overline{\boldsymbol{M}}_{k}^{2}=\overline{\boldsymbol{X}}\left(i_{k}^{\prime}, 2\right) & , k=1, \ldots, p \\
\overline{\boldsymbol{N}}_{\boldsymbol{k}}^{1}=\overline{\boldsymbol{X}}\left(i_{k}^{\prime \prime}, 1\right), \overline{\boldsymbol{N}}_{\boldsymbol{k}}^{2}=\overline{\boldsymbol{X}}\left(i_{k}^{\prime \prime}, 2\right) & , k=1, \ldots, q
\end{aligned}
$$

where $\overline{\boldsymbol{X}}\left(i_{k}^{\prime}, 1\right), \overline{\boldsymbol{X}}\left(i_{k}^{\prime}, 2\right), \overline{\boldsymbol{X}}\left(i_{k}^{\prime \prime}, 1\right)$ and $\overline{\boldsymbol{X}}\left(i_{k}^{\prime \prime}, 2\right)$ are given by Definitions 5 and 6 . Moreover, let $\underline{N}=\left(\underline{n}_{i j}\right)_{q \times n}$ be a $q \times n$ matrix whose entries are defined as follows:

$$
\underline{n}_{k j}=\left\{\begin{array}{cc}
b_{i_{k}^{\prime \prime}} & j \in J_{i_{k}^{\prime \prime}} \cup\left\{i_{k}^{\prime \prime}\right\} \\
-\infty & \text { otherwise }
\end{array} \quad, k=1, \ldots, q ~~,~~ j=1, \ldots, n\right.
$$

\textbf{Remark 2.} For each $i_{k}^{\prime} \in I_{2}(k=1, \ldots, p)$, the $k^{\prime}$ th row of matrix $\bar{M}^{1}\left(\bar{M}^{2}\right)$ is exactly the maximal solution $\overline{\boldsymbol{X}}\left(i_{k}^{\prime}, 1\right)\left(\overline{\boldsymbol{X}}\left(i_{k}^{\prime}, 2\right)\right)$ of $S\left(\boldsymbol{a}_{i_{k}^{\prime}}, b_{i_{k}}\right)$. Similarly, for each $i_{k}^{\prime \prime} \in I_{3}(k=1, \ldots, q)$, the $k^{\prime}$ th row of matrix $\bar{N}^{1}\left(\bar{N}^{2}\right)$ is exactly the maximal solution $\bar{X}\left(i_{k}^{\prime \prime}, 1\right)\left(\bar{X}\left(i_{k}^{\prime \prime}, 2\right)\right)$ of $S\left(a_{i_{k}^{\prime}}, b_{i_{k}^{*}}\right)$. Moreover, Matrix $\underline{N}$ basically plays the same role as do matrix $M^{1^{*}}$ defined in [14], the minimal solution matrix $\breve{\Gamma}$ introduced in [44] and the simplified matrix presented in [31]. We refer the reader to [14] in which some details were provided about the relationships between the three mentioned matrices.\\

\textbf{Lemma 5.} (Rule 1). Let $I_{2}=\left\{i_{1}^{\prime}, \ldots, i_{p}^{\prime}\right\}$ and consider a fixed $i_{k}^{\prime} \in I_{2}$. Also, suppose that there exists some $j \in J$ such that $\max \left\{\left(\underline{\boldsymbol{X}}_{1}\right)_{j},\left(\underline{\boldsymbol{X}}_{2}\right)_{j}\right\}>\bar{m}_{k j}^{1}\left(\max \left\{\left(\underline{\boldsymbol{X}}_{1}\right)_{j},\left(\underline{\boldsymbol{X}}_{2}\right)_{j}\right\}>\bar{m}_{k j}^{2}\right)$, where $\left(\underline{X}_{1}\right)_{j}$ and $\left(\underline{X}_{2}\right)_{j}$ denote the $j$ 'th components of $\underline{X}_{1}$ and $\underline{\boldsymbol{X}}_{2}$, respectively. Then, each triple $\left(\underline{e}, e^{\prime}, e^{\prime \prime}\right)$ such that $\vec{e}^{\prime}\left(i_{k}^{\prime}\right)=1\left(\vec{e}^{\prime}\left(i_{k}^{\prime}\right)=2\right)$ is inadmissible.\\

\textbf{Proof.} Let $\max \left\{\left(\underline{\boldsymbol{X}}_{1}\right)_{j},\left(\underline{\boldsymbol{X}}_{2}\right)_{j}\right\}>\bar{m}_{k j}^{1}$ for some $j \in J$. Also, consider a triple $\left(\underline{e}, e^{\prime}, e^{\prime \prime}\right)$ in which $\underline{e}$ and $e^{\prime \prime}$ are two arbitrary functions in $\underline{E}$ and $E^{\prime \prime}$, respectively, and $e^{\prime} \in E^{\prime}$ such that $e^{\prime}\left(i_{k}^{\prime}\right)=1$. From Definition 9 , we have $\overline{\boldsymbol{X}}_{2}\left(e^{\prime}\right)=\min _{i \in I_{2}}\left\{\overline{\boldsymbol{X}}\left(i, e^{\prime}(i)\right)\right\}=\min \left\{\overline{\boldsymbol{X}}\left(i_{k}^{\prime}, 1\right), \min _{i \in I_{2}-\left\{i_{k}^{\prime}\right\}}\left\{\overline{\boldsymbol{X}}\left(i, e^{\prime}(i)\right)\right\}\right\}$. Therefore, $j^{\prime}$ th component of $\overline{\boldsymbol{X}}_{2}\left(e^{\prime}\right)$ is obtained as 
\begin{align*}
\overline{\boldsymbol{X}}_{2}\left(e^{\prime}\right)_{j} &= \min \left\{\overline{\boldsymbol{X}}\left(i_{k}^{\prime}, 1\right)_{j}, \min _{i \in I_{2}-\left\{i_{k}^{\prime}\right\}}\left\{\overline{\boldsymbol{X}}\left(i, e^{\prime}(i)\right)_{j}\right\}\right\},
\end{align*}
that implies $\overline{\boldsymbol{X}}_{2}\left(e^{\prime}\right)_{j} \leq \overline{\boldsymbol{X}}\left(i_{k}^{\prime}, 1\right)_{j}$. Since $\overline{\boldsymbol{X}}\left(i_{k}^{\prime}, 1\right)_{j}=\bar{m}_{k j}^{1}$ (from Definition 11) and $\bar{m}_{k j}^{1}<\max \left\{\left(\underline{\boldsymbol{X}}_{1}\right)_{j},\left(\underline{\boldsymbol{X}}_{2}\right)_{j}\right\}$, we have $\quad \overline{\boldsymbol{X}}_{2}\left(e^{\prime}\right)_{j}<\max \left\{\left(\underline{\boldsymbol{X}}_{1}\right)_{j},\left(\underline{\boldsymbol{X}}_{2}\right)_{j}\right\}$. Thus, $\max \left\{\underline{\boldsymbol{X}}_{1}, \underline{\boldsymbol{X}}_{2}\right\} \nleq \overline{\boldsymbol{X}}_{2}\left(e^{\prime}\right)$ that implies $\max \left\{\underline{\boldsymbol{X}}_{1}, \underline{\boldsymbol{X}}_{2}, \underline{\boldsymbol{X}}_{3}(\underline{e})\right\} \nleq \overline{\boldsymbol{X}}_{2}\left(e^{\prime}\right)$. Consequently, $\max \left\{\underline{\boldsymbol{X}}_{1}, \underline{\boldsymbol{X}}_{2}, \underline{\boldsymbol{X}}_{3}(\underline{e})\right\} \nleq \min \left\{\overline{\boldsymbol{X}}_{1}, \overline{\boldsymbol{X}}_{2}\left(e^{\prime}\right), \overline{\boldsymbol{X}}_{3}\left(e^{\prime \prime}\right)\right\}$ that violates Relation (5). Hence, by Definition 10, it follows that $\left(\underline{e}, e^{\prime}, e^{\prime \prime}\right) \notin T$. The proof is similar if $\max \left\{\left(\underline{\boldsymbol{X}}_{1}\right)_{j},\left(\underline{\boldsymbol{X}}_{2}\right)_{j}\right\}>\bar{m}_{k j}^{2}$ for some $j \in J$ and $e^{\prime}\left(i_{k}^{\prime}\right)=2$. \\

\textbf{Remark 3.} By Definition 8 and Remark $1, e^{\prime}\left(i_{k}^{\prime}\right) \in \operatorname{dom}\left(e^{\prime}\left(i_{k}^{\prime}\right)\right)=\{1,2\}$ for each $e^{\prime} \in E^{\prime}$ and each $i_{k}^{\prime} \in I_{2}$. In the case that $e^{\prime}\left(i_{k}^{\prime}\right)=1$, we have $\overline{\boldsymbol{X}}\left(i_{k}^{\prime}, e^{\prime}\left(i_{k}^{\prime}\right)\right)=\overline{\boldsymbol{X}}\left(i_{k}^{\prime}, 1\right)$ and therefore according to Definition 9 , the maximal solution $\overline{\boldsymbol{X}}\left(i_{k}^{\prime}, 1\right)$ participates in the construction of $\overline{\boldsymbol{X}}_{2}\left(e^{\prime}\right)=\min _{i \in I_{2}}\left\{\overline{\boldsymbol{X}}\left(i, e^{\prime}(i)\right)\right\}$. However, if $\max \left\{\left(\underline{\boldsymbol{X}}_{1}\right)_{j},\left(\underline{\boldsymbol{X}}_{2}\right)_{j}\right\}>\bar{m}_{k j}^{1}$, Lemma 5 implies that each triple $\left(\underline{e}, e^{\prime}, e^{\prime \prime}\right)$ such that $e^{\prime}\left(i_{k}^{\prime}\right)=1$ is inadmissible and it leads to an empty cell
\begin{align*}
\left[\max \left\{\underline{\boldsymbol{X}}_{1}, \underline{\boldsymbol{X}}_{2}, \underline{\boldsymbol{X}}_{3}(e)\right\}, \min \left\{\overline{\boldsymbol{X}}_{1}, \overline{\boldsymbol{X}}_{2}\left(e^{\prime}\right), \overline{\boldsymbol{X}}_{3}\left(e^{\prime \prime}\right)\right\}\right]. 
\end{align*}
Therefore, based on Rule 1 , if $\max \left\{\left(\underline{\boldsymbol{X}}_{1}\right)_{j},\left(\underline{\boldsymbol{X}}_{2}\right)_{j}\right\}>\bar{m}_{k j}^{1}$ for some $i_{k}^{\prime} \in I_{2}$ and $j \in J$, we set $\overline{\boldsymbol{M}}_{\boldsymbol{k}}^{1}=[\infty, \infty, \ldots, \infty]_{1 \times n}$ and $\operatorname{dom}\left(e^{\prime}\left(i_{k}^{\prime}\right)\right)=\{2\}$, that is, the domain of $e^{\prime}\left(i_{k}^{\prime}\right)$ is decreased from $\{1,2\}$ to $\{2\}$. The replacement $\overline{\boldsymbol{M}}_{k}^{1}=[\infty, \infty, \ldots, \infty]_{1 \times n}$ means that the maximal solution $\overline{\boldsymbol{X}}\left(i_{k}^{\prime}, 1\right)$ (i.e., $\overline{\boldsymbol{X}}\left(i_{k}^{\prime}, e^{\prime}\left(i_{k}^{\prime}\right)\right)$ with $\left.e^{\prime}\left(i_{k}^{\prime}\right)=1\right)$ is not allowed to be selected in the construction of $\overline{\boldsymbol{X}}_{2}\left(e^{\prime}\right)$. Similar notation is used if $\max \left\{\left(\underline{\boldsymbol{X}}_{1}\right)_{j},\left(\underline{\boldsymbol{X}}_{2}\right)_{j}\right\}>\bar{m}_{k j}^{2}$ for some $i_{k}^{\prime} \in I_{2}$ and $j \in J$. In this case, we set $\overline{\boldsymbol{M}}_{\boldsymbol{k}}^{2}=[\infty, \infty, \ldots, \infty]_{1 \times n}$ and $\operatorname{dom}\left(e^{\prime}\left(i_{k}^{\prime}\right)\right)=\{1\}$.\\

\textbf{Corollary 4.} If $\overline{\boldsymbol{M}}_{k}^{1}=\overline{\boldsymbol{M}}_{k}^{2}=[\infty, \infty, \ldots, \infty]_{1 \times n}$ for some $k \in\{1, \ldots, p\}$, then Problem (2) is infeasible.\\
\textbf{Proof.} From Remark 3, the assignments $\overline{\boldsymbol{M}}_{\boldsymbol{k}}^{1}=\overline{\boldsymbol{M}}_{\boldsymbol{k}}^{2}=[\infty, \infty, \ldots, \infty]_{1 \times n}$ mean that each value of $e^{\prime}\left(i_{k}^{\prime}\right)$, whether $e^{\prime}\left(i_{k}^{\prime}\right)=1$ or $e^{\prime}\left(i_{k}^{\prime}\right)=2$, leads to an inadmissible triple. Therefore, each triple $\left(\underline{e}, e^{\prime}, e^{\prime \prime}\right)$ violates Relation (5). Now, Theorem 3 implies $S(\boldsymbol{A}, \boldsymbol{b}) \neq \varnothing$. \\

\textbf{Lemma 6.} (Rule 2). Let $I_{3}=\left\{i_{1}^{\prime \prime}, \ldots, i_{q}^{\prime \prime}\right\}$ and consider a fixed $i_{k}^{\prime \prime} \in I_{3}$. Also, suppose that there exists some $j \in J$ such that $\max \left\{\left(\underline{\boldsymbol{X}}_{1}\right)_{j},\left(\underline{\boldsymbol{X}}_{2}\right)_{j}\right\}>\bar{n}_{k j}^{1}\left(\max \left\{\left(\underline{\boldsymbol{X}}_{1}\right)_{j},\left(\underline{\boldsymbol{X}}_{2}\right)_{j}\right\}>\bar{n}_{k j}^{2}\right)$. Then, each triple $\left(\underline{e}, e^{\prime}, e^{\prime \prime}\right)$ such that $e^{\prime \prime}\left(i_{k}^{\prime \prime}\right)=1\left(e^{\prime \prime}\left(i_{k}^{\prime \prime}\right)=2\right)$ is inadmissible.\\
\textbf{Proof.} The proof is similar to the proof of Lemma 5 by replacing $I_{2}, E^{\prime}, e^{\prime}$ and $\overline{\boldsymbol{X}}_{2}\left(e^{\prime}\right)$ with $I_{3}, E^{\prime \prime}, e^{\prime \prime}$ and $\overline{\boldsymbol{X}}_{3}\left(e^{\prime \prime}\right)$, respectively.\\

\textbf{Remark 4.} According to Rule 2, if $\max \left\{\left(\underline{\boldsymbol{X}}_{1}\right)_{j},\left(\underline{\boldsymbol{X}}_{2}\right)_{j}\right\}>\bar{n}_{k j}^{1}$ for some $i_{k}^{\prime \prime} \in I_{3}$ and $j \in J$, we set $\overline{\boldsymbol{N}}_{\boldsymbol{k}}^{1}=[\infty, \infty, \ldots, \infty]_{1 \times n}$ and $\operatorname{dom}\left(e^{\prime \prime}\left(i_{k}^{\prime \prime}\right)\right)=\{2\}$. Similarly, if $\max \left\{\left(\underline{\boldsymbol{X}}_{1}\right)_{j},\left(\underline{\boldsymbol{X}}_{2}\right)_{j}\right\}>\bar{n}_{k j}^{2}$ for some $i_{k}^{\prime \prime} \in I_{3}$ and $j \in J$, we set $\bar{N}_{k}^{2}=[\infty, \infty, \ldots, \infty]_{1 \times n}$ and $\operatorname{dom}\left(e^{\prime \prime}\left(i_{k}^{\prime \prime}\right)\right)=\{1\}$.\\

\textbf{Corollary 5.} If $\overline{\boldsymbol{N}}_{k}^{1}=\overline{\boldsymbol{N}}_{k}^{2}=[\infty, \infty, \ldots, \infty]_{1 \times n}$ for some $k \in\{1, \ldots, q\}$, then Problem (2) is infeasible.\\
\textbf{Proof.} The proof is similar to the proof of Corollary 4 by replacing $I_{2}, E^{\prime}$ and $e^{\prime}$ with $I_{3}, E^{\prime \prime}$ and $e^{\prime \prime}$, respectively.\\

\textbf{Lemma 7} (Rule 3). Let $I_{3}=\left\{i_{1}^{\prime \prime}, \ldots, i_{q}^{\prime \prime}\right\}$ and consider a fixed $i_{k}^{\prime \prime} \in I_{3}$. Also, suppose that there exists some $j \in J_{i_{k}^{\prime \prime}}$ such that $\underline{n}_{k j}>\left(\overline{\boldsymbol{X}}_{1}\right)_{j}$, where $\left(\overline{\boldsymbol{X}}_{1}\right)_{j}$ denotes the $j$ 'th component of $\bar{X}_{1}$. Then, each triple $\left(\underline{e}, e^{\prime}, e^{\prime \prime}\right)$ such that $\underline{e}\left(i_{k}^{\prime \prime}\right)=j$ is inadmissible.\\
\textbf{Proof.} Suppose that $\left(\underline{e}, e^{\prime}, e^{\prime \prime}\right)$ in which $e^{\prime}$ and $e^{\prime \prime}$ are two arbitrary functions in $E^{\prime}$ and $E^{\prime \prime}$, respectively, and $\underline{e} \in \underline{E}$ such that $\underline{e}\left(i_{k}^{\prime \prime}\right)=j$. From Definition 9, we have $\underline{\boldsymbol{X}}_{3}(\underline{e})=\max _{i \in I_{3}}\{\underline{\boldsymbol{X}}(i, \underline{e}(i))\}=\max \left\{\underline{\boldsymbol{X}}\left(i_{k}^{\prime \prime}, j\right), \max _{i \in I_{3}-\left\{i_{k}^{\prime}\right\}}\{\underline{\boldsymbol{X}}(i, \underline{e}(i))\}\right\}$. Therefore, $j$ 'th component of $\underline{\boldsymbol{X}}_{3}(\underline{e})$ is obtained as
\begin{align*}
\underline{\boldsymbol{X}}_{3}(\underline{e})_{j}=\max \left\{\underline{\boldsymbol{X}}\left(i_{k}^{\prime \prime}, j\right)_{j}, \max_{i \in I_{3} - {\left\{k_{k}^{\prime \prime}\right\}}}\left\{\underline{\boldsymbol{X}}(i, \underline{e}(i))_{j}\right\}\right\},
\end{align*}
that implies $\underline{\boldsymbol{X}}_{3}(\underline{e})_{j} \geq \underline{\boldsymbol{X}}\left(i_{k}^{\prime \prime}, j\right)_{j}$. On the other hand, since $\underline{\boldsymbol{X}}\left(i_{k}^{\prime \prime}, j\right)_{j}=\underline{n}_{k j}$ (from Definitions 6 and 11) and $\underline{n}_{k j}>\left(\overline{\boldsymbol{X}}_{1}\right)_{j}$, it follows that $\underline{\boldsymbol{X}}_{3}(\underline{e})_{j}>\left(\overline{\boldsymbol{X}}_{1}\right)_{j}$. Thus, $\underline{\boldsymbol{X}}_{3}(\underline{e}) \nleq \overline{\boldsymbol{X}}_{1}$, and therefore $\max \left\{\underline{\boldsymbol{X}}_{1}, \underline{\boldsymbol{X}}_{2}, \underline{\boldsymbol{X}}_{3}(\underline{e})\right\} \nleq \overline{\boldsymbol{X}}_{1}$. Hence, we have $\max \left\{\underline{\boldsymbol{X}}_{1}, \underline{\boldsymbol{X}}_{2}, \underline{\boldsymbol{X}}_{3}(\underline{e})\right\} \nleq \min \left\{\overline{\boldsymbol{X}}_{1}, \overline{\boldsymbol{X}}_{2}\left(e^{\prime}\right), \overline{\boldsymbol{X}}_{3}\left(e^{\prime \prime}\right)\right\}$ that violates Relation (5). Now, from Definition 10, we have $\left(\underline{e}, e^{\prime}, e^{\prime \prime}\right) \notin T$ that completes the proof.\\

\textbf{Remark 5.} Based on Rule 3, if $\underline{n}_{k j}>\left(\overline{\boldsymbol{X}}_{1}\right)_{j}$ for some $i_{k}^{\prime \prime} \in I_{3}$ and $j \in J_{i_{k}^{\prime \prime}}$, then we set $\underline{n}_{k j}=-\infty$ and $\operatorname{dom}\left(\underline{e}\left(i_{k}^{\prime \prime}\right)\right)=J_{i_{k}^{\prime \prime}}-\{j\}$. The assignment $\underline{n}_{k j}=-\infty$ means that the minimal solution $\underline{\boldsymbol{X}}\left(i_{k}^{\prime \prime}, j\right)$ (i.e., $\underline{\boldsymbol{X}}\left(i_{k}^{\prime \prime}, \underline{e}\left(i_{k}^{\prime \prime}\right)\right)$ with $\underline{e}\left(i_{k}^{\prime \prime}\right)=j$ ) is not allowed to be selected in the construction of $\underline{\boldsymbol{X}}_{3}(\underline{e})=\max _{i \in I_{3}}\{\underline{\boldsymbol{X}}(i, \underline{e}(i))\}$ as defined in Definition 9.\\

\textbf{Corollary 6.} Suppose that $i_{k}^{\prime \prime} \in I_{3}$ and $\underline{n}_{k j}=-\infty, \forall j \in J_{i_{k}^{\prime \prime}}$. Then, Problem (2) is infeasible.\\
\textbf{Proof.} By Definition 8 and Remark 1, $\operatorname{dom}\left(\underline{e}\left(i_{k}^{\prime \prime}\right)\right)=J_{i_{k}^{\prime \prime}}$. Moreover, since $\underline{n}_{k j}=-\infty$, $\forall j \in J_{i_{k}^{\prime \prime}}$, Lemma 7 implies that any possible value of $\underline{e}\left(i_{k}^{\prime \prime}\right)$ leads to an inadmissible triple $\left(\underline{e}, e^{\prime}, e^{\prime \prime}\right)$ that violates Relation (5). Thus, from Theorem 3, we have $S(\boldsymbol{A}, \boldsymbol{b})=\varnothing$.\\

\textbf{Lemma 8.} (Rule 4). Let $I_{2}=\left\{i_{1}^{\prime}, \ldots, i_{p}^{\prime}\right\}$ and $I_{3}=\left\{i_{1}^{\prime \prime}, \ldots, i_{q}^{\prime \prime}\right\}$. Also, suppose that $i_{r}^{\prime} \in I_{2}$ and $i_{s}^{\prime \prime} \in I_{3}$ such that $a_{i_{r}^{\prime} r_{s}^{\prime \prime}}>b_{i_{r}^{\prime}}$ and $b_{i_{r}^{\prime \prime}}<b_{i_{s}^{\prime \prime}}$. Then, each triple $\left(\underline{e}, e^{\prime}, e^{\prime \prime}\right)$ such that $e^{\prime}\left(i_{r}^{\prime}\right)=2$ is inadmissible.\\
\textbf{Proof.} Consider a triple $\left(\underline{e}, e^{\prime}, e^{\prime \prime}\right)$ in which $\underline{e}$ and $e^{\prime \prime}$ are two arbitrary functions in $\underline{E}$ and $E^{\prime \prime}$, respectively, and $e^{\prime} \in E^{\prime}$ such that $e^{\prime}\left(i_{r}^{\prime}\right)=2$. Similar to the proof of Lemmas 5 and $7, i_{s}^{\prime \prime}$ th components of $\overline{\boldsymbol{X}}_{2}\left(e^{\prime}\right)$ and $\underline{\boldsymbol{X}}_{3}(\underline{e})$ are obtained as follows:
$$
\begin{aligned}
& \overline{\boldsymbol{X}}_{2}\left(e^{\prime}\right)_{i_{s}^{\prime}}=\min _{i \in I_{2}}\left\{\overline{\boldsymbol{X}}\left(i, e^{\prime}(i)\right)_{i_{s}^{\prime \prime}}\right\}=\min \left\{\overline{\boldsymbol{X}}\left(i_{r}^{\prime}, 2\right)_{i_{i^{\prime}}}, \min _{i \in I_{2}-\{{i^{\prime}_{r}}\}}\left\{\overline{\boldsymbol{X}}\left(i, e^{\prime}(i)\right)_{i_{s}^{\prime \prime}}\right\}\right\} \\
& \underline{\boldsymbol{X}}_{3}(\underline{e})_{i_{s}^{\prime \prime}}=\max _{i \in I_{3}}\left\{\underline{\boldsymbol{X}}(i, \underline{e}(i))_{i_{s}^{\prime \prime}}\right\}=\max \left\{\underline{\boldsymbol{X}}\left(i_{s}^{\prime \prime}, j\right)_{i_{s}^{\prime \prime}}, \max _{i \in I_{3}-\{i^{\prime \prime}_{s}\}}\left\{\underline{\boldsymbol{X}}(i, \underline{e}(i))_{i_{s}^{\prime \prime}}\right\}\right\}.
\end{aligned}
$$
Therefore, we have
\begin{equation}
    \begin{array}{l}
        \overline{\boldsymbol{X}}_{2}\left(e^{\prime}\right)_{i_{s}^{\prime \prime}} \leq \overline{\boldsymbol{X}}\left(i_{r}^{\prime}, 2\right)_{i_{s}^{\prime \prime}} 
    \end{array}
\end{equation}
\begin{equation}
    \begin{array}{l}
        \underline{\boldsymbol{X}}_{3}(\underline{e})_{i_{s}^{\prime \prime}} \geq \underline{\boldsymbol{X}}\left(i_{s}^{\prime \prime}, j\right)_{i_{s}^{\prime \prime}}.
    \end{array}
\end{equation}
Also, we have $\overline{\boldsymbol{X}}\left(i_{r}^{\prime}, 2\right)_{i_{s}^{\prime \prime}}=\bar{m}_{r r_{s}^{\prime \prime}}^{2}=b_{i_{r}^{\prime}}$ (Definitions 5 and 11) and $\underline{\boldsymbol{X}}\left(i_{s}^{\prime \prime}, j\right)_{i_{i^{\prime \prime}}}=\underline{n}_{s_{i}^{\prime \prime}}=b_{i_{s}^{\prime \prime}}$ (Definitions 6 and 11). The last equalities together with (6) and (7) imply $\overline{\boldsymbol{X}}_{2}\left(e^{\prime}\right)_{i_{s}^{\prime \prime}} \leq b_{i_{t}^{\prime}}$ and $\underline{\boldsymbol{X}}_{3}(\underline{e})_{i_{s}^{\prime}} \geq b_{i_{s}^{\prime \prime}}$. Therefore, $\quad \overline{\boldsymbol{X}}_{2}\left(e^{\prime}\right)_{i_{s}^{\prime \prime}} \leq b_{i_{i}^{\prime}}<b_{i_{s}^{\prime \prime}} \leq \underline{\boldsymbol{X}}_{3}(\underline{e})_{i_{s}^{\prime \prime}}$, and then $\underline{\boldsymbol{X}}_{3}(\underline{e}) \nleq \overline{\boldsymbol{X}}_{2}\left(e^{\prime}\right)$. So, $\max \left\{\underline{\boldsymbol{X}}_{1}, \underline{\boldsymbol{X}}_{2}, \underline{\boldsymbol{X}}_{3}(\underline{e})\right\} \nleq \min \left\{\overline{\boldsymbol{X}}_{1}, \overline{\boldsymbol{X}}_{2}\left(e^{\prime}\right), \overline{\boldsymbol{X}}_{3}\left(e^{\prime \prime}\right)\right\}$ that violates Relation (5). Now, from Definition 10, we conclude $\left(\underline{e}, e^{\prime}, e^{\prime \prime}\right) \notin T$, and the result is obtained.\\

\textbf{Lemma 9.} (Rule 5). Let $I_{3}=\left\{i_{1}^{\prime \prime}, \ldots, i_{q}^{\prime \prime}\right\}$. Also, suppose that $i_{r}^{\prime \prime}, i_{s}^{\prime \prime} \in I_{3}$ such that $i_{r}^{\prime \prime} \neq i_{s}^{\prime \prime}$, $a_{i_{r}^{\prime \prime} i_{s}^{\prime \prime}}>b_{i_{r}^{\prime \prime}}$ and $b_{i_{r}^{\prime \prime}}<b_{i_{s}^{\prime \prime}}$. Then, each triple $\left(\underline{e}, e^{\prime}, e^{\prime \prime}\right)$ such that $e^{\prime \prime}\left(i_{r}^{\prime \prime}\right)=2$ is inadmissible.\\
\textbf{Proof.} The proof is quite similar to the proof of Lemma 8.\\

\textbf{Remark 6.} According to Rule 4, if there exist $i_{r}^{\prime} \in I_{2}$ and $i_{s}^{\prime \prime} \in I_{3}$ such that $a_{i^{\prime} r_{s}^{\prime \prime}}>b_{i_{r}^{\prime}}$ and $b_{i_{r}^{\prime}}<b_{i_{s}^{\prime \prime}}$, then we set $\overline{\boldsymbol{M}}_{r}^{2}=[\infty, \infty, \ldots, \infty]_{1 \times n}$ and $\operatorname{dom}\left(e^{\prime}\left(i_{r}^{\prime}\right)\right)=\{1\}$. Similarly, based on Rule 5, if there exist $i_{r}^{\prime \prime}, i_{s}^{\prime \prime} \in I_{3}$ such that $i_{r}^{\prime \prime} \neq i_{s}^{\prime \prime}, \quad a_{i_{r}^{\prime \prime} s_{s}^{\prime \prime}}>b_{i_{r}^{\prime \prime}}$ and $b_{i_{r}^{\prime \prime}}<b_{i_{s}^{\prime \prime}}$, then we set $\bar{N}_{r}^{2}=[\infty, \infty, \ldots, \infty]_{1 \times n}$ and $\operatorname{dom}\left(e^{\prime \prime}\left(i_{r}^{\prime \prime}\right)\right)=\{1\}$.\\

\textbf{Lemma 10.} (Rule 6). Suppose that Rules $1-5$ have been already applied. Also, suppose that $\operatorname{dom}\left(e^{\prime}\left(i_{r}^{\prime}\right)\right)=\{1\}$ for some $i_{r}^{\prime} \in I_{2}$ (or equivalently, $\overline{\boldsymbol{M}}_{r}^{1} \neq[\infty, \infty, \ldots, \infty]_{1 \times n}$ and $\overline{\boldsymbol{M}}_{r}^{2}=[\infty, \infty, \ldots, \infty]_{1 \times n}$). If there exists $i_{s}^{\prime \prime} \in I_{3}$ such that $i_{r}^{\prime} \in J_{i_{s}^{\prime \prime}}$ and $b_{i_{r}^{\prime}}<b_{i_{s}^{\prime \prime}}$, Then, each triple $\left(\underline{e}, e^{\prime}, e^{\prime \prime}\right)$ such that $\underline{e}\left(i_{s}^{\prime \prime}\right)=i_{r}^{\prime}$ is inadmissible.\\
\textbf{Proof.} Suppose that $\left(\underline{e}, e^{\prime}, e^{\prime \prime}\right)$ in which $e^{\prime}$ and $e^{\prime \prime}$ are two arbitrary functions in $E^{\prime}$ and $E^{\prime \prime}$, respectively, and $\underline{e} \in \underline{E}$ such that $\underline{e}\left(i_{s}^{\prime \prime}\right)=i_{r}^{\prime}$. Since $\operatorname{dom}\left(e^{\prime}\left(i_{r}^{\prime}\right)\right)=\{1\}$, then there exists only one option for the value of $e^{\prime}\left(i_{r}^{\prime}\right)$, i.e., $e^{\prime}\left(i_{r}^{\prime}\right)=1$. The $i_{r}^{\prime}$ th components of $\overline{\boldsymbol{X}}_{2}\left(e^{\prime}\right)$ and $\underline{\boldsymbol{X}}_{3}(\underline{e})$ are obtained as follows:
$$
\begin{aligned}
& \overline{\boldsymbol{X}}_{2}\left(e^{\prime}\right)_{i_{r^{\prime}}^{\prime}}=\min _{i \in I_{2}}\left\{\overline{\boldsymbol{X}}\left(i, e^{\prime}(i)\right)_{i_{r^{\prime}}}\right\}=\min \left\{\overline{\boldsymbol{X}}\left(i_{r}^{\prime}, 1\right)_{i^{\prime}_{r}}, \min _{i \in I_{2}-\left\{i_{r^{\prime}}\right\}}\left\{\bar{X}\left(i, e^{\prime}(i)\right)_{i_{i_{r}^{\prime}}}\right\}\right\} \\
& \underline{\boldsymbol{X}}_{3}(\underline{e})_{i_{r}^{\prime}}=\max _{i \in I_{3}}\left\{\underline{\boldsymbol{X}}(i, \underline{e}(i))_{i_{i_{r}^{\prime}}}\right\}=\max \left\{\underline{\boldsymbol{X}}\left(i_{s}^{\prime \prime}, i_{r}^{\prime}\right)_{i^{\prime}_{r}}, \max _{i \in I_{3}-\left\{i_{s}^{\prime \prime}\right\}}\left\{\underline{\boldsymbol{X}}(i, \underline{e}(i))_{i_{r}^{\prime}}\right\}\right\}
\end{aligned}
$$
Therefore, we have
\begin{equation}
    \begin{array}{l}
        \overline{\boldsymbol{X}}_{2}\left(e^{\prime}\right)_{i_{r}^{\prime}} \leq \overline{\boldsymbol{X}}\left(i_{r}^{\prime}, 1\right)_{i_{r}^{\prime}}
    \end{array}
\end{equation}
\begin{equation}
    \begin{array}{l}
        \underline{\boldsymbol{X}}_{3}(\underline{e})_{i_{r}^{\prime}} \ge \max\underline{\boldsymbol{X}}\left(i_{s}^{\prime \prime}, i_{r}^{\prime}\right)_{i^{\prime}_{r}}
    \end{array}
\end{equation}
Also, we have $\overline{\boldsymbol{X}}\left(i_{r}^{\prime}, 1\right)_{i_{r}^{\prime}}=\bar{m}_{r i_{r}^{\prime}}^{1}=b_{i_{r}^{\prime}}$ (Definitions 5 and 11) and $\underline{\boldsymbol{X}}\left(i_{s}^{\prime \prime}, i_{r}^{\prime}\right)_{i_{r}^{\prime}}=\underline{n}_{s i_{r}^{\prime}}=b_{i_{s}^{\prime \prime}}$ (Definitions 6 and 11). The latter equalities together with (8) and (9) imply $\overline{\boldsymbol{X}}_{2}\left(e^{\prime}\right)_{i_{r}^{\prime}} \leq b_{i_{r}^{\prime}}$ and $\underline{\boldsymbol{X}}_{3}(\underline{e})_{i_{r}^{\prime}} \geq b_{i_{s}^{\prime \prime}}$. Therefore, $\overline{\boldsymbol{X}}_{2}\left(e^{\prime}\right)_{i_{r}^{\prime}} \leq b_{i_{r}^{\prime}}<b_{i_{s}^{\prime \prime}} \leq \underline{\boldsymbol{X}}_{3}(\underline{e})_{i_{r}^{\prime}}$, and then $\underline{\boldsymbol{X}}_{3}(\underline{e}) \nleq \overline{\boldsymbol{X}}_{2}\left(e^{\prime}\right)$. So, $\max \left\{\underline{\boldsymbol{X}}_{1}, \underline{\boldsymbol{X}}_{2}, \underline{\boldsymbol{X}}_{3}(\underline{e})\right\} \nleq \min \left\{\overline{\boldsymbol{X}}_{1}, \overline{\boldsymbol{X}}_{2}\left(e^{\prime}\right), \overline{\boldsymbol{X}}_{3}\left(e^{\prime \prime}\right)\right\}$ that violates Relation (5). Hence, we have $\left(\underline{e}, e^{\prime}, e^{\prime \prime}\right) \notin T$ from Definition 10 that completes the proof.\\

\textbf{Lemma 11.} (Rule 7). Suppose that Rules $1-5$ have been already applied. Also, suppose that $\operatorname{dom}\left(e^{\prime \prime}\left(i_{r}^{\prime \prime}\right)\right)=\{1\}$ for some $i_{r}^{\prime \prime} \in I_{3}$ (or equivalently, $\bar{N}_{r}^{1} \neq[\infty, \infty, \ldots, \infty]_{1 \times n}$ and $\left.\bar{N}_{r}^{2}=[\infty, \infty, \ldots, \infty]_{1 \times n}\right)$. If there exists $i_{s}^{\prime \prime} \in I_{3}$ such that $i_{r}^{\prime \prime} \neq i_{s}^{\prime \prime}, i_{r}^{\prime \prime} \in J_{i_{s}^{\prime \prime}}$ and $b_{i_{r}^{\prime \prime}}<b_{i_{s}^{\prime \prime}}$, then each triple $\left(\underline{e}, e^{\prime}, e^{\prime \prime}\right)$ such that $\underline{e}\left(i_{s}^{\prime \prime}\right)=i_{r}^{\prime \prime}$ is inadmissible.\\
\textbf{Proof.} The proof is quite similar to the proof of Lemma 10.\\

\textbf{Remark 7.} According to Rule 6 , if there exist $i_{r}^{\prime} \in I_{2}$ and $i_{s}^{\prime \prime} \in I_{3}$ such that $\operatorname{dom}\left(e^{\prime}\left(i_{r}^{\prime}\right)\right)=\{1\}, \quad i_{r}^{\prime} \in J_{i_{s}^{\prime \prime}}$ and $b_{i_{r}^{\prime}}<b_{i_{s}^{\prime \prime}}$, we set $\underline{n}_{s i_{r}^{\prime}}=-\infty$ and $\operatorname{dom}\left(\underline{e}\left(i_{s}^{\prime \prime}\right)\right)=J_{i_{s}^{\prime \prime}}-\left\{i_{r}^{\prime}\right\}$. Similarly, based on Rule 7 , if there exist $i_{r}^{\prime \prime}, i_{s}^{\prime \prime} \in I_{3}$ such that $i_{r}^{\prime \prime} \neq i_{s}^{\prime \prime}, \operatorname{dom}\left(e^{\prime \prime}\left(i_{r}^{\prime \prime}\right)\right)=\{1\}$, $i_{r}^{\prime \prime} \in J_{i_{s}^{\prime \prime}}$ and $b_{i_{r}^{\prime \prime}}<b_{i_{s}^{\prime \prime}}$, we set $\underline{n}_{s i_{r}^{\prime \prime}}=-\infty$ and $\operatorname{dom}\left(\underline{e}\left(i_{s}^{\prime \prime}\right)\right)=J_{i_{s}^{\prime \prime}}-\left\{i_{r}^{\prime \prime}\right\}$.\\

\section{Optimal solution of the problem}
This section describes the optimal solutions of Problem (2). Also, it is proved that under some assumptions on the components of matrix $\boldsymbol{A}$ and vector $\boldsymbol{b}$, Problem (2) automatically yields binary optimal solutions.

\textbf{Definition 12.} Let $J^{+}=\left\{j \in J: c_{j} \geq 0\right\}$ and $J^{-}=\left\{j \in J: c_{j}<0\right\}$. Associated with each triple $e=\left(\underline{e}, e^{\prime}, e^{\prime \prime}\right)$ in $T$, we define a vector $\boldsymbol{x}_{e}=\left[\left(x_{e}\right)_{1},\left(x_{e}\right)_{2}, \ldots,\left(x_{e}\right)_{n}\right]$ as follows:

\begin{equation}
    \begin{array}{l}
    \left(x_{e}\right)_{j}= \begin{cases}\max \left\{\left(\underline{\boldsymbol{X}}_{1}\right)_{j},\left(\underline{\boldsymbol{X}}_{2}\right)_{j}, \underline{\boldsymbol{X}}_{3}(\underline{e}_{j})\right\} & , j \in J^{+} \\ 
    \min \left\{\left(\overline{\boldsymbol{X}}_{1}\right)_{j}, \overline{\boldsymbol{X}}_{2}\left(e^{\prime}\right)_{j}, \overline{\boldsymbol{X}}_{3}\left(e^{\prime \prime}\right)_{j}\right\}
    & , j \in J^{-}\end{cases}
    \end{array}
\end{equation}

According to the following theorem, the optimal solution of Problem (2) is always attained as $\boldsymbol{x}_{e}$ for some $e=\left(\underline{e}, e^{\prime}, e^{\prime \prime}\right) \in T$.\\

\textbf{Theorem 4.} Let $\boldsymbol{c}^{T} \boldsymbol{x}_{e^{*}}=\min \left\{\boldsymbol{c}^{T} \boldsymbol{x}_{e}: e \in T\right\}$. Then, $\boldsymbol{x}_{e^{*}}$ is the optimal solution of Problem (2).\\
\textbf{Proof.} Let $\boldsymbol{x} \in S(\boldsymbol{A}, \boldsymbol{b})$. Therefore, Theorem 2 implies that
\begin{align*}
\boldsymbol{x} \in\left[\max \left\{\underline{\boldsymbol{X}}_{1}, \underline{\boldsymbol{X}}_{2}, \underline{\boldsymbol{X}}_{3}(\underline{e})\right\}, \min \left\{\overline{\boldsymbol{X}}_{1}, \overline{\boldsymbol{X}}_{2}\left(e^{\prime}\right), \overline{\boldsymbol{X}}_{3}\left(e^{\prime \prime}\right)\right\}\right]
\end{align*}
for some $e=\left(\underline{e}, e^{\prime}, e^{\prime \prime}\right) \in T$. Therefore, from Definition 12, $\left(x_{e}\right)_{j} \leq x_{j}, \quad \forall j \in J^{+}$, and $x_{j} \leq\left(x_{e}\right)_{j}, \forall j \in J^{-}$. Hence, $\sum_{j \in J^{+}} c_{j}\left(x_{e}\right)_{j} \leq \sum_{j \in J^{+}} c_{j^{\prime}} x_{j} $ and $ \sum_{j \in J^{-}} c_{j}\left(x_{e}\right)_{j} \leq \sum_{j \in J^{-}} c_{j} x_{j}$, and therefore $\boldsymbol{c}^{T} \boldsymbol{x}_{e}=\sum_{j \in J^{+}} c_{j}\left(x_{e}\right)_{j}+\sum_{j \in J^{-}} c_{j}\left(x_{e}\right)_{j} \leq \sum_{j \in J^{+}} c_{j} x_{j}+\sum_{j \in J^{-}} c_{j} x_{j}=\boldsymbol{c}^{T} \boldsymbol{x}$. Now, by the assumption of the theorem, it follows that $\boldsymbol{c}^{T} \boldsymbol{x}_{e^{*}} \leq \boldsymbol{c}^{T} \boldsymbol{x}_{e} \leq \boldsymbol{c}^{T} \boldsymbol{x}$.\\

\textbf{Corollary 7.} If Problem (2) is expressed as a maximization problem, the global optimal solution $\boldsymbol{x}^{*}$ is obtained by $\boldsymbol{c}^{T} \boldsymbol{x}^{*}=\max \left\{\boldsymbol{c}^{T} \boldsymbol{x}_{e}: e \in T\right\}$ where solutions $\boldsymbol{x}_{e}=\left[\left(x_{e}\right)_{1},\left(x_{e}\right)_{2}, \ldots,\left(x_{e}\right)_{n}\right]\left(e=\left(\underline{e}, e^{\prime}, e^{\prime \prime}\right) \in T\right)$ are defined as follows
\begin{equation}
    \begin{array}{l}
        \left(x_{e}\right)_{j}= \begin{cases}
        \min \left\{\left(\overline{\boldsymbol{X}}_{1}\right)_{j}, \overline{\boldsymbol{X}}_{2}\left(e^{\prime}\right)_{j}, \overline{\boldsymbol{X}}_{3}\left(e^{\prime \prime}\right)_{j}\right\} & , j \in J^{+} \\ 
        
        \max \left\{\left(\underline{\boldsymbol{X}}_{1}\right)_{j}, \left(\underline{\boldsymbol{X}}_{2}\right)_{j}, \underline{\boldsymbol{X}}_{3}(\underline{e})_{j}\right\} & , j \in J^{-}
        \end{cases}
    \end{array}
\end{equation}
The following theorem shows that in the special cases of Problem (2), optimal solutions are binary valued.\\

\textbf{Theorem 5.} Consider Problem (2) where $\boldsymbol{A}=\left(a_{i j}\right)_{m \times n}$ is a matrix such that $a_{i j} \in\{0,1\} \quad(\forall i \in I$ and $\forall j \in J)$ and $\boldsymbol{b}=\mathbf{0}_{m \times 1}$ is a zero vector. If $\boldsymbol{x}^{*}=\left(x_{1}^{*}, x_{2}^{*}, \ldots, x_{n}^{*}\right)$ is the optimal solution of the problem, then $x_{j}^{*} \in\{0,1\}, \forall j \in J$.\\
\textbf{Proof.} By contradiction, suppose that $0<x_{j_{0}}^{*}<1$ for some $j_{0} \in J$. Since $\boldsymbol{b}=\mathbf{0}$, then from Relation (3), the constraints of the problem are expressed as $\max _{j=1}^{n}\left\{\min \left\{a_{i j}, x_{i}, x_{j}\right\}\right\}=0, \forall i \in I$. Therefore, $\min \left\{a_{i j}, x_{i}, x_{j}\right\}=0, \forall i \in I \quad$ and $\forall j \in J$. Particularly, we have
\begin{equation}
    \begin{array}{l}
        \min \left\{a_{i j_{0}}, x_{i}^{*}, x_{j_{0}}^{*}\right\}=0 , \forall i \in I \\
        \min \left\{a_{j_{0} j}, x_{j_{0}}^{*}, x_{j}^{*}\right\}=0 , \forall j \in J \text { such that } j_{0} \leq j
    \end{array}
\end{equation}
The assumption $0<x_{j_{0}}^{*}<1$ and (12) result in the following two conditions:
\begin{equation}
    \begin{array}{l}
        (I) ~~ \text {for each} ~~ i \in I, ~~ \text {either} ~~ a_{i j_{0}}=0 ~~ \text {or} ~~ x_{i}^{*}=0.\\
        (II) ~~ \text {for each} ~~ j \in J ~~ \text {such that} ~~ j_{0} \leq j, ~~\text{either}~~ a_{j_{j} j}=0 ~~\text{or}~~ x_{j}^{*}=0.
    \end{array}
\end{equation}
At the first case, consider $j_{0} \in J^{-}$and define solution $x^{\prime}$ such that $x_{j_{0}}^{\prime}=1$ and $x_{j}^{\prime}=x_{j}^{*}$ , $\forall j \in J-\left\{j_{0}\right\}$. Based on (13), it is clear that $\min \left\{a_{i j_{0}}, x_{i}^{\prime}, x_{j_{0}}^{\prime}\right\}=\min \left\{a_{i j_{0}}, x_{i}^{*}, 1\right\}=0, \forall i \in I$, and $\min \left\{a_{j_{0} j}, x_{j_{0}}^{\prime}, x_{j}^{\prime}\right\}=\min \left\{a_{j_{0} j}, 1, x_{j}^{*}\right\}=0, \forall j \in J$ such that $j_{0} \leq j$; that is, $x^{\prime}$ satisfies all the equations stated in (12). Moreover, since $x_{j}^{\prime}=x_{j}^{*}, \forall j \in J-\left\{j_{0}\right\}, \quad x^{\prime}$ also satisfies the other equations in which $x_{j_{0}}^{*}$ does not appear. So, $\boldsymbol{x}^{\prime}$ is a feasible solution to the problem. On the other hand, we have $c_{j_{0}} x_{j_{0}}^{\prime}<c_{j_{0}} x_{j_{0}}^{*}$ and $\sum_{j \in J-\left\{j_{0}\right\}} c_{j} x_{j}^{\prime}=\sum_{j \in J-\left\{j_{0}\right\}} c_{j} x_{j}^{*}$. Therefore, $\boldsymbol{c}^{T} \boldsymbol{x}^{\prime}=c_{j_{0}} x_{j_{0}}^{\prime}+\sum_{j \in J-\left\{j_{0}\right\}} c_{j} x_{j}^{\prime}<c_{j_{0}} x_{j_{0}}^{*}+\sum_{j \in J-\left\{j_{0}\right\}} c_{j} x_{j}^{*}=\boldsymbol{c}^{T} \boldsymbol{x}^{*}$ that violates the optimality of $\boldsymbol{x}^{*}$. If $j_{0} \in J^{+}$, the proof is simpler by defining solution $\boldsymbol{x}^{\prime}$ such that $x_{j_{0}}^{\prime}=0$ and $x_{j}^{\prime}=x_{j}^{*}$ ,$\forall j \in J-\left\{j_{0}\right\}$.\\

\textbf{Corollary 8.} Suppose that Problem (2) is expressed as a maximization problem where $\boldsymbol{A}=\left(a_{i j}\right)_{m \times n}$ is a matrix such that $a_{i j} \in\{0,1\} \quad(\forall i \in I$ and $\forall j \in J)$ and $\boldsymbol{b}=\mathbf{0}_{m \times 1}$ is a zero vector. If $\boldsymbol{x}^{*}=\left(x_{1}^{*}, x_{2}^{*}, \ldots, x_{n}^{*}\right)$ is the optimal solution of the problem, then $x_{j}^{*} \in\{0,1\}, \forall j \in J$.\\
The following algorithm summarizes the preceding discussion.\\

\textbf{Algorithm 1.} (Solution of Problem (2)).\\
Given Problem (2) and suppose that $I_{2}=\left\{i_{1}^{\prime}, \ldots, i_{p}^{\prime}\right\}$ and $I_{3}=\left\{i_{1}^{\prime \prime}, \ldots, i_{q}^{\prime \prime}\right\}$. \\

1. If $J_{i}=\varnothing$ for some $i \in I$, then stop; Problem (2) is infeasible (Corollary 1 ).
\begin{enumerate}
  \setcounter{enumi}{1}
  \item Compute solutions $\underline{X}_{1}, \bar{X}_{1}$ and $\underline{X}_{2}$ (Definition 9).

  \item If $\max \left\{\underline{\boldsymbol{X}}_{1}, \underline{\boldsymbol{X}}_{2}\right\} \nless \overline{\boldsymbol{X}}_{1}$, then stop; Problem (2) is infeasible (Theorem 1).

  \item Apply Rule 1 in order to reduce matrices $\overline{\boldsymbol{M}}^{1}$ and $\overline{\boldsymbol{M}}^{2}$. If $\overline{\boldsymbol{M}}_{\boldsymbol{k}}^{1}=\overline{\boldsymbol{M}}_{\boldsymbol{k}}^{2}=[\infty, \infty, \ldots, \infty]$ for some $k \in\{1, \ldots, p\}$, then stop; Problem (2) is infeasible (Corollary 4).

  \item Apply Rule 2 in order to reduce matrices $\bar{N}^{1}$ and $\bar{N}^{2}$. If $\bar{N}_{k}^{1}=\bar{N}_{k}^{2}=[\infty, \infty, \ldots, \infty]$ for some $k \in\{1, \ldots, q\}$, then stop; problem (2) is infeasible (Corollary 5).

  \item Apply Rule 3 in order to reduce matrix $\underline{N}$. If there exists some $k \in\{1, \ldots, q\}$ such that $\underline{n}_{k j}=-\infty, \forall j \in J_{i_{k}^{\prime \prime}}$, then stop; problem (2) is infeasible (Corollary 6).

  \item Apply Rules 4 and 5 in order to reduce matrices $\bar{M}^{2}$ and $\bar{N}^{2}$.

  \item Apply Rules 6 and 7 in order to reduce matrix $\underline{N}$.

  \item For each triple $\left(\underline{e}, e^{\prime}, e^{\prime \prime}\right)$ selected based on the reduced matrices $\bar{M}^{1}, \bar{M}^{2}, \bar{N}^{1}$, $\overline{\boldsymbol{N}}^{2}$ and $\underline{\boldsymbol{N}}$, compute solutions $\overline{\boldsymbol{X}}_{2}\left(e^{\prime}\right), \overline{\boldsymbol{X}}_{3}\left(e^{\prime \prime}\right)$ and $\underline{\boldsymbol{X}}_{3}(\underline{e})$ (Remarks 2-7).

  \item If $\underline{\boldsymbol{X}}_{3}(\underline{e}), \overline{\boldsymbol{X}}_{2}\left(e^{\prime}\right)$ and $\overline{\boldsymbol{X}}_{3}\left(e^{\prime \prime}\right)$ satisfy Relation (5), generate solution $\boldsymbol{x}_{e}$.

  \item Find the optimal solution $\boldsymbol{x}_{e^{*}}$ by $\boldsymbol{c}^{T} \boldsymbol{x}_{e^{*}}=\min \left\{\boldsymbol{c}^{T} \boldsymbol{x}_{e}: e \in T\right\}$ (Theorem 4).

\end{enumerate}

\section{An important especial case of Problem (2); Minimum Vertex Cover Problem}
Consider an undirected simple connected graph $G=(V, E)$ consisting of a set $V$ of $n$ vertices (nodes) and a set $E$ of edges whose elements are unordered pairs of the distinct vertices. Let $(i, j)$ denote an undirected edge between two vertices $i$ and $j$. Formally, a vertex cover $V^{\prime}$ of an undirected graph $G=(V, E)$ is a subset of $V$ such that $(i, j) \in E$ implies $i \in V^{\prime}$ and $j \in V^{\prime}$; that is, every edge has at least one endpoint in $V^{\prime}$. Such a set is said to cover the edges of $G$. Figure 3 shows two examples of vertex covers, where vertices of $V^{\prime}$ have been marked in black. A minimal vertex cover is a vertex cover with the smallest possible size. Figure $3(b)$ shows an example of minimum vertex covers.\\

\begin{figure}[ht]
    \begin{center}
	    \includegraphics[height=5cm]{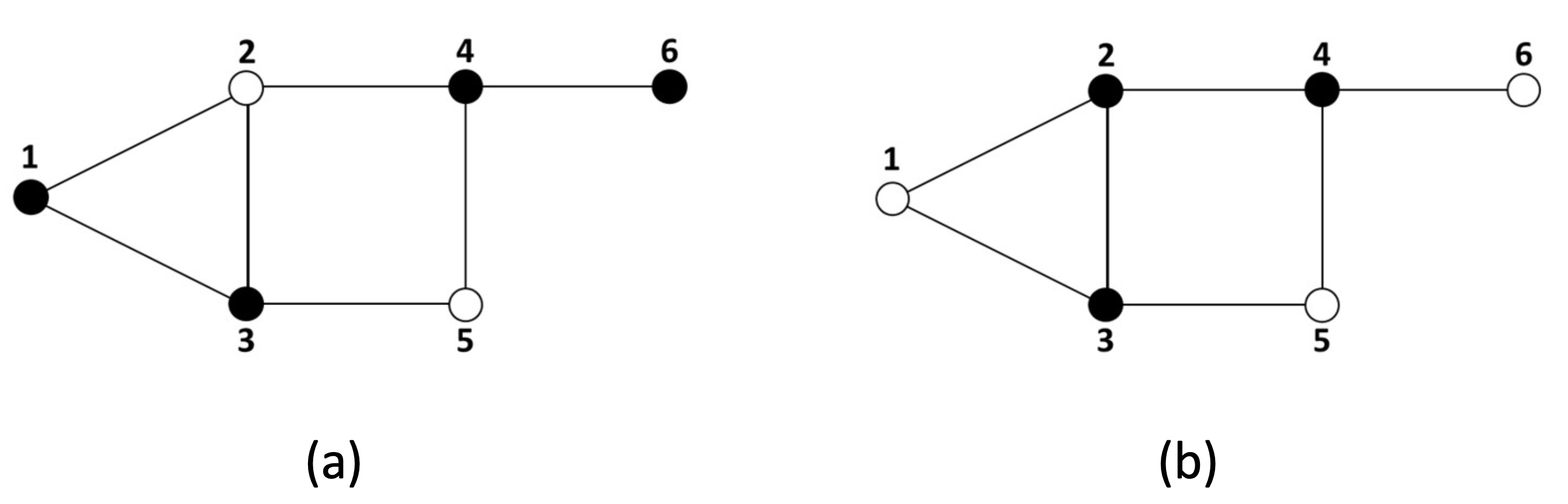}
	    \caption{(a) Example of a vertex cover, (b) Example of a minimum vertex cover.}
	\end{center}
\end{figure} 

\subsection{Formulation of the minimum vertex cover problem}
Let the node-node adjacency matrix $\boldsymbol{A}=\left(a_{i j}\right)_{n \times n}$ have a row and a column corresponding to every vertex, and its $(i, j)$ 'th entry equals 1 (i.e., $a_{i j}=1$ ) if $(i, j) \in E$ and equals 0 (i.e., $\left.a_{i j}=0\right)$ if $(i, j) \notin E$. Particularly, we have $a_{i i}=0$ for each $i \in I$. Moreover, associate with each vertex $j \in V$ a binary variable $v_{j}^{\prime} \in\{0,1\}$ such that $v_{j}^{\prime}=1$ if $j \in V^{\prime}$ (black vertices in Fig. 3) and $v_{j}^{\prime}=0$ otherwise (white vertices). Based on the above statements, it is clear that $V^{\prime}$ is a vertex cover if and only if there are no vertices $i$ and $j$ such that $a_{i j}=1$ and $v_{i}^{\prime}=v_{j}^{\prime}=0$. So, the minimal vertex cover problem can be formulated as follows:
\begin{equation}
    \begin{array}{l}
         \min Z=\sum_{j \in J} v_{j}^{\prime} \\
        \max _{j=1}^{n}\left\{\min \left\{a_{i j}, 1-v_{i}^{\prime}, 1-  v_{j}^{\prime}\right\}\right\}=0 \quad, \quad i \in I, j \in J \\
         v_{j}^{\prime} \in\{0,1\} \quad, \quad j \in J
    \end{array}
\end{equation}
where $a_{i j} \in\{0,1\}(\forall i \in I$ and $\forall j \in J)$. It is worth noting that if $a_{i j}=1$ for any vertices $i$ and $j$, the above equations prevent the variables $v_{i}^{\prime}$ and $v_{j}^{\prime}$ from having zero value at the same time. Therefore, each feasible solution of (14) corresponds to a vertex cover $V^{\prime}$, and vice versa (vertex $j$ belongs to $V^{\prime}$ iff $v_{j}^{\prime}=1$ ). So, the goal of the above problem is to find a solution with the maximum number of variables $v_{j}^{\prime}$ with zero values (or equivalently, the minimum number of variables $v_{j}^{\prime}$ with one values). By setting $x_{j}=1-v_{j}^{\prime}$, Problem (14) is converted into an equivalent problem as follows:
\begin{equation}
    \begin{array}{l}
         \min Z=\sum_{j \in J}\left(1-x_{j}\right) \\
         \max _{j=1}^{n}\left\{\min \left\{a_{i j}, x_{i}, x_{j}\right\}\right\}=0 \quad, \quad i \in I, j \in J \\
         x_{j} \in\{0,1\} \quad, \quad j \in J
    \end{array}
\end{equation}
Now, consider the following problem derived from (15) by manipulating the objective function and replacing the constraints $x_{j} \in\{0,1\}$ by $x_{j} \in[0,1](j \in J)$ :
\begin{equation}
    \begin{array}{l}
            \max Z  =\sum_{j \in J} x_{j} \\
            \max _{j=1}^{n}\left\{\min \left\{a_{i j}, x_{i}, x_{j}\right\}\right\}=0 \quad, \quad i \in I, j \in J  \\
            \boldsymbol{x} \in[0,1]^{n} 
    \end{array}
\end{equation}

Clearly, Problem (16) is a special case of Problem (2) (see Corollaries 7 and 8) in which $b_{i}=0(i \in I)$ and $c_{j}=1 (j \in J)$. Hence, according to Corollary 8, all optimal solutions of Problem (16) are binary, and therefore problems (15) and (16) have the same optimal solutions. As a consequence, the minimum vertex cover problem is indeed a special case of Problem (2).

\subsection{Properties of Problem (16)}
As mentioned before, $G=(V, E)$ is a simple graph (i.e., a graph without loops and multiple edges). So, in Problem (16) we have $a_{i i}=0$ for each $i \in I$, and therefore it is concluded from Definition 3 that $I_{1}=I_{3}=\varnothing$ and $I=I_{2}$. Consequently, $E^{\prime \prime}=\underline{E}=\varnothing$ (see Definition 8) and $S\left(a_{i}, b_{i}\right)=\varnothing, \forall i \in I_{1} \cup I_{2}$, which in turn implies that $S_{1}(\boldsymbol{A}, \boldsymbol{b})=S_{3}(\boldsymbol{A}, \boldsymbol{b})=\varnothing$ (see Definition 7). Moreover, since $I=I_{2}$ and $S(\boldsymbol{A}, \boldsymbol{b})=\bigcap_{i \in I} S\left(\boldsymbol{a}_{i}, b_{i}\right)$, from Definition 7 we have
\begin{equation}
    \begin{array}{l}
        S_{2}(\boldsymbol{A}, \boldsymbol{b})=\bigcap_{i \in I_{2}} S\left(\boldsymbol{a}_{i}, b_{i}\right)=\bigcap_{i \in I} S\left(\boldsymbol{a}_{i}, b_{i}\right)=S(\boldsymbol{A}, \boldsymbol{b})
    \end{array}
\end{equation}
which together with Lemma 4(b) imply
\begin{equation}
    \begin{array}{l}
        S(\boldsymbol{A}, \boldsymbol{b})=S_{2}(\boldsymbol{A}, \boldsymbol{b})=\bigcup_{e^{\prime} \in E^{\prime}}\left[\underline{\boldsymbol{X}}_{2}, \overline{\boldsymbol{X}}_{2}\left(e^{\prime}\right)\right]
    \end{array}
\end{equation}
According to (18), $S(\boldsymbol{A}, \boldsymbol{b}) \neq \varnothing$ iff $\left[\underline{\boldsymbol{X}}_{2}, \overline{\boldsymbol{X}}_{2}\left(e^{\prime}\right)\right]$ for some $e^{\prime} \in E^{\prime}$. Hence, function $e^{\prime} \in E^{\prime}$ is admissible, if $\left[\underline{\boldsymbol{X}}_{2}, \overline{\boldsymbol{X}}_{2}\left(e^{\prime}\right)\right] \neq \varnothing$ or equivalently $\underline{\boldsymbol{X}}_{2} \leq \overline{\boldsymbol{X}}_{2}\left(e^{\prime}\right)$. By noting this fact and the equalities $E^{\prime \prime}=\underline{E}=\varnothing,(5)$ is reduced to the following relation:
\begin{equation}
    \begin{array}{l}
        \underline{\boldsymbol{X}}_{2} \leq \overline{\boldsymbol{X}}_{2}\left(e^{\prime}\right)
    \end{array}
\end{equation}
and set $T$ (defined in Definition 10) is also modified as follows:
\begin{equation}
    \begin{array}{l}
        T=\left\{e^{\prime} \in E^{\prime}: \underline{\boldsymbol{X}}_{2} \leq \overline{\boldsymbol{X}}_{2}\left(e^{\prime}\right)\right\}
    \end{array}
\end{equation}
Additionally, from Definition 5, the following results are directly obtained for each $i \in I_{2}$ and $j \in J$ :
\begin{equation}
    \begin{array}{l}
        \bar{X}(i, 1)_{j}=\left\{\begin{array}{ll}
        0 & j=i \\
        1 & j \neq i
        \end{array} \quad, \quad \bar{X}(i, 2)_{j}=\left\{\begin{array}{ll}
        0 & a_{i j}=1 \\
        1 & a_{i j}=0
\end{array} \quad, \quad \underline{X}(i)_{j}=0\right.\right.
    \end{array}
\end{equation}
Furthermore, the equalities $I_{1}=I_{3}=\varnothing$ and Definition $11 \mathrm{imply} \overline{\boldsymbol{N}}^{1}=\overline{\boldsymbol{N}}^{2}=\underline{\boldsymbol{N}}=\varnothing$. Also, from (21) and Definition 11, we have $\overline{\boldsymbol{M}}^{\mathbf{1}}=\mathbf{1}-\mathbf{I}$ where $\mathbf{1}_{n \times n}$ is a matrix of ones and $\mathbf{I}_{n \times n}$ is the identity matrix, and
\begin{equation}
    \begin{array}{l}
        \overline{\boldsymbol{M}}^{2}=\left[\begin{array}{ccc}
        1-a_{11} & \cdots & 1-a_{1 n} \\
        \vdots & \ddots & \vdots \\
        1-a_{n 1} & \cdots & 1-a_{n n}
        \end{array}\right]
    \end{array}
\end{equation}
For the sake of expository reference, some especial results have been summarized in the following corollary:\\

\textbf{Corollary 9.} Consider Problem (16) and let $i \in I_{2}=I$.\\
\textbf{(a)} $\overline{\boldsymbol{X}}(i, 1)$ contains a zero in the $i^{\prime}$ th position and ones everywhere else. \textbf{(b)} $\overline{\boldsymbol{X}}(i, 2)^{T}=\left(1-a_{i 1}, 1-a_{i 2}, \ldots, 1-a_{i n}\right)$. \textbf{(c)} $\underline{\boldsymbol{X}}(i)=\mathbf{0}$ where $\mathbf{0}$ denotes the zero vector. \textbf{(d)} $\overline{\boldsymbol{X}}\left(i, e^{\prime}(i)\right)_{j} \in\{0,1\}, \forall e^{\prime} \in E^{\prime}$ and $\forall j \in J$. \textbf{(e)} $\underline{\boldsymbol{X}}_{2}=\mathbf{0}$. \textbf{(f)} $T=E^{\prime}$. \textbf{(g)} $S(\boldsymbol{A}, \boldsymbol{b})=\bigcup_{e^{\prime} \in E^{\prime}}\left[\mathbf{0}, \overline{\boldsymbol{X}}_{2}\left(e^{\prime}\right)\right]$.\\
\textbf{Proof.} \textbf{(a)-(c)} The results follow from (21). \textbf{(d)} Since $e^{\prime}(i) \in\{1,2\}$ (see Definition 8), the proof is directly resulted from Parts (a) and (b). \textbf{(e)} The result is easily attained from Definition 9 and (21). \textbf{(f)} From Part (e) and (20), we have $T=\left\{e^{\prime} \in E^{\prime}: \mathbf{0} \leq \overline{\boldsymbol{X}}_{2}\left(e^{\prime}\right)\right\}$. However, since the condition $\mathbf{0} \leq \overline{\boldsymbol{X}}_{2}\left(e^{\prime}\right)$ is satisfied for each $e^{\prime} \in E^{\prime}$, then $T=E^{\prime}$. \textbf{(g)} The result follows from Part (e) and (18). \\ 
Since $T=E^{\prime}$ (from Corollary 9(f)), then for each $e \in T$ we can consider $e=e^{\prime}$ where $e^{\prime} \in E^{\prime}$. On the other hand, for Problem (16), it is easily found that $J^{-}=\varnothing$ and $J^{+}=J$ (because, in the objective function we have $c_{j}=1, \forall j \in J$ ). Consequently, by setting $\boldsymbol{c}=\mathbf{1}$ where $\mathbf{1}_{n \times 1}$ denotes the sum vector (i.e., a vector having each component equal to one), Corollary 7 is reduced to Corollary 10 below.\\

\textbf{Corollary 10.} Consider Problem (16) and let $\mathbf{1}_{n \times 1}$ denote the sum vector. Then, the global optimal solution $\boldsymbol{x}^{*}$ is obtained by $\mathbf{1}^{T} \boldsymbol{x}^{*}=\max \left\{\mathbf{1}^{T} \boldsymbol{x}_{e^{\prime}}: e^{\prime} \in E^{\prime}\right\}$ where solutions $\boldsymbol{x}_{e^{\prime}}=\left[\left(x_{e^{\prime}}\right)_{1},\left(x_{e^{\prime}}\right)_{2}, \ldots,\left(x_{e^{\prime}}\right)_{n}\right]\left(e^{\prime} \in E^{\prime}\right)$ are defined as follows:
\begin{equation}
    \begin{array}{l}    \left(x_{e^{\prime}}\right)_{j}=\overline{\boldsymbol{X}}_{2}\left(e^{\prime}\right)_{j} \quad, j \in J
    \end{array}
\end{equation}

\textbf{Lemma 12.} Let $x^{*}$ denote the global optimal solution of Problem (16). Then, (a). $\mathbf{1}^{T} \boldsymbol{x}^{*}=\max \left\{\mathbf{1}^{T} \overline{\boldsymbol{X}}_{2}\left(e^{\prime}\right): e^{\prime} \in E^{\prime}\right\}$. (b) $\boldsymbol{x}^{*}=\overline{\boldsymbol{X}}_{2}\left(e^{*}\right)$ for some $e^{*} \in E^{\prime}$. (c). $x_{j}^{*} \in\{0,1\}$, $\forall j \in J$.\\
\textbf{Proof.} \textbf{(a)} From Corollary 10, we have $\mathbf{1}^{T} \boldsymbol{x}^{*}=\max \left\{\mathbf{1}^{T} \boldsymbol{x}_{e^{\prime}}: e^{\prime} \in E^{\prime}\right\}$. Now, the result follows from (23). \textbf{(b)} It is a direct consequence of Part (a). \textbf{(c)} The proof is directly resulted from Part (b), Definition 9 and Corollary 9(d).\\

\textbf{Lemma 13.} Suppose that $\boldsymbol{x}^{*}=\overline{\boldsymbol{X}}_{2}\left(e^{*}\right)\left(e^{*} \in E^{\prime}\right)$ is the global optimal solution of Problem (16) and $e^{\prime} \in E^{\prime}$.

\textbf{(a)}. If $e^{\prime}(i)=1, \forall i \in I_{2}$, then $\overline{\boldsymbol{X}}_{2}\left(e^{\prime}\right)=\mathbf{0}$.

\textbf{(b)}. Let $e^{\prime}\left(i_{0}\right)=2$ for some $i_{0} \in I_{2}$ and $a_{i_{0} j}=1$ (or equivalently, $\bar{m}_{i_{0} j}^{2}=0$ ). Then, $\overline{\boldsymbol{X}}_{2}\left(e^{\prime}\right)_{j}=0$.

\textbf{(c)}. Let $e^{*}\left(i_{0}\right)=2$ for some $i_{0} \in I_{2}$. Then, $e^{*}(i)=1$ for each $i \in I_{2}-\left\{i_{0}\right\}$ such that $a_{i_{0} i}=1$ (or equivalently, $\bar{m}_{i_{0} i}^{2}=0$ ). \\
\textbf{Proof.} \textbf{(a)} It is sufficient to show that $\overline{\boldsymbol{X}}_{2}\left(e^{\prime}\right)_{j}=0, \forall j \in J$. From Definition 9, $\overline{\boldsymbol{X}}_{2}\left(e^{\prime}\right)=\min _{i \in I_{2}}\left\{\overline{\boldsymbol{X}}\left(i, e^{\prime}(i)\right)\right\}$. Since $I_{2}=I$ and $e^{\prime}(i)=1, \forall i \in I_{2}$, then $\overline{\boldsymbol{X}}_{2}\left(e^{\prime}\right)=\min _{i \in I}\{\overline{\boldsymbol{X}}(i, 1)\}$. Thus, for each $j \in J, \overline{\boldsymbol{X}}_{2}\left(e^{\prime}\right)_{j}=\min _{i \in I}\left\{\overline{\boldsymbol{X}}(i, 1)_{j}\right\}$. Now, Corollary $9(\mathrm{a})$ implies that $\overline{\boldsymbol{X}}(j, 1)_{j}=0$, and therefore $\quad \overline{\boldsymbol{X}}_{2}\left(e^{\prime}\right)_{j}=\min _{i \in I}\left\{\overline{\boldsymbol{X}}(i, 1)_{j}\right\}=\min \left\{\overline{\boldsymbol{X}}(j, 1)_{j}, \min _{i \in I-\{j\}}\left\{\overline{\boldsymbol{X}}(i, 1)_{j}\right\}\right\}=\min \left\{0, \min _{i \in I-\{j\}}\left\{\overline{\boldsymbol{X}}(i, 1)_{j}\right\}\right\}=0$. \textbf{(b)} Similar to Part (a), from Definition 9 and the equality $I_{2}=I$, we have
\begin{align*}
    \overline{\boldsymbol{X}}_{2}\left(e^{\prime}\right)_{j}=\min _{i \in I_{2}}\left\{\overline{\boldsymbol{X}}\left(i, e^{\prime}(i)\right)_{j}\right\}=\min _{i \in I}\left\{\overline{\boldsymbol{X}}\left(i, e^{\prime}(i)\right)_{j}\right\}
    \\ =\min \left\{\overline{\boldsymbol{X}}\left(i_{0}, 2\right)_{j}, \min _{i \in I-\left\{i_{0}\right\}}\left\{\overline{\boldsymbol{X}}\left(i, e^{\prime}(i)\right)_{j}\right\}\right\}.
\end{align*}
But, Corollary 9(b) requires that $\overline{\boldsymbol{X}}\left(i_{0}, 2\right)_{j}=1-a_{i_{0} j}=0$. Therefore, $\overline{\boldsymbol{X}}_{2}\left(e^{\prime}\right)_{j}=\min \left\{0, \min _{i \in I- 
\{i_{0}\}}\left\{\overline{\boldsymbol{X}}\left(i, e^{\prime}(i)\right)_{j}\right\}\right\}=0$. \textbf{(c)} Let $e^{*}\left(i_{0}\right)=2$ and $e^{*}(i)=1$ for each $i \in I_{2}-\left\{i_{0}\right\}$ such that $a_{i_{0}}=1$. Also, suppose that $i^{\prime} \in I_{2}-\left\{i_{0}\right\}$ and $a_{i_{0} i^{\prime}}=1$. Now, consider $e^{\prime} \in E^{\prime}$ such that $e^{\prime}\left(i^{\prime}\right)=2$ and $e^{\prime}(i)=e^{*}(i)$, $\forall i \in I_{2}-\left\{i^{\prime}\right\}$. We show that $\mathbf{1}^{T} \overline{\boldsymbol{X}}_{2}\left(e^{*}\right) \geq \mathbf{1}^{T} \overline{\boldsymbol{X}}_{2}\left(e^{\prime}\right)$; that is, $\overline{\boldsymbol{X}}_{2}\left(e^{\prime}\right)$ cannot be a better solution than $\overline{\boldsymbol{X}}_{2}\left(e^{*}\right)$. Firstly, we note that since $e^{\prime}\left(i_{0}\right)=e^{*}\left(i_{0}\right)=2$ and $a_{i_{0} i^{\prime}}=1$, we have from Part (b) that $\overline{\boldsymbol{X}}_{2}\left(e^{\prime}\right)_{i^{\prime}}=\overline{\boldsymbol{X}}_{2}\left(e^{*}\right)_{i^{\prime}}=0$. Otherwise, assume that $j \neq i^{\prime}$. So, since $\overline{\boldsymbol{X}}\left(i^{\prime}, 2\right)_{j}=1-a_{i^{\prime}}$ (Corollary 9(b)), then $\overline{\boldsymbol{X}}\left(i^{\prime}, 2\right)_{j}=1$ if $a_{i^{\prime} j}=0$, and $\overline{\boldsymbol{X}}\left(i^{\prime}, 2\right)_{j}=0$ if $a_{i^{\prime} j}=1$. Therefore, we have\\ 
$$
\begin{aligned}
\overline{\boldsymbol{X}}_{2}\left(e^{*}\right)_{j} & =\min _{i \in I}\left\{\overline{\boldsymbol{X}}\left(i, e^{*}(i)\right)_{j}\right\}=\min \left\{\overline{\boldsymbol{X}}\left(i^{\prime}, e^{*}\left(i^{\prime}\right)\right)_{j}, \min _{i \in I-\left\{i^{\prime}\right\}}\left\{\overline{\boldsymbol{X}}\left(i, e^{*}(i)\right)_{j}\right\}\right\} \\
& =\min \left\{\overline{\boldsymbol{X}}\left(i^{\prime}, 1\right)_{j}, \min _{i \in I-\{i^{\prime}\}}\left\{\overline{\boldsymbol{X}}\left(i, e^{*}(i)\right)_{j}\right\}\right\}=\min \left\{1, \min _{i \in I-\{i^{\prime}\}}\left\{\overline{\boldsymbol{X}}\left(i, e^{*}(i)\right)_{j}\right\}\right\} \\
& =\min _{i \in I-\left\{i^{\prime}\right\}}\left\{\overline{\boldsymbol{X}}\left(i, e^{*}(i)\right)_{j}\right\}=\min _{i \in I-\left\{i^{\prime}\right\}}\left\{\overline{\boldsymbol{X}}\left(i, e^{\prime}(i)\right)_{j}\right\} \\ 
& \geq \min \left\{\overline{\boldsymbol{X}}\left(i^{\prime}, e^{\prime}\left(i^{\prime}\right)\right)_{j}, \min _{i \in I-\left\{i^{\prime}\right\}}\left\{\overline{\boldsymbol{X}}\left(i, e^{\prime}(i)\right)_{j}\right\}\right\} \\
& \geq \min \left\{\overline{\boldsymbol{X}}\left(i^{\prime}, 2\right)_{j}, \min _{i \in I-\left\{i^{\prime}\right\}}\left\{\overline{\boldsymbol{X}}\left(i, e^{\prime}(i)\right)_{j}\right\}\right\} \\
& =\overline{\boldsymbol{X}}_{2}\left(e^{\prime}\right)_{j}
\end{aligned}
$$
Therefore, $\overline{\boldsymbol{X}}_{2}\left(e^{*}\right)_{j} \geq \overline{\boldsymbol{X}}_{2}\left(e^{\prime}\right)_{j}, \forall j \in J$, that implies $\mathbf{1}^{T} \overline{\boldsymbol{X}}_{2}\left(e^{*}\right) \geq \mathbf{1}^{T} \overline{\boldsymbol{X}}_{2}\left(e^{\prime}\right)$.\\

\textbf{Theorem 6.} Suppose that $\boldsymbol{x}^{*}=\overline{\boldsymbol{X}}_{2}\left(e^{*}\right)\left(e^{*} \in \bar{E}_{2}\right)$ is the global optimal solution of Problem (16).

\textbf{(a)}. There exists at least one $i_{0} \in I_{2}$ such that $e^{*}\left(i_{0}\right)=2$.

\textbf{(b)}. If $e^{*}\left(i_{1}\right)=e^{*}\left(i_{2}\right)=2$, then $a_{i_{1} i_{2}}=a_{i_{2} i_{1}}=0$ (or equivalently, $\bar{m}_{i_{i} i_{2}}^{2}=\bar{m}_{i_{2} i_{1}}^{2}=1$ ). \\
\textbf{Proof.} \textbf{(a)} By contradiction, suppose that $e^{*}(i)=1, \forall i \in I_{2}$. Hence, Lemma 13(a) implies that $\boldsymbol{x}^{*}=\overline{\boldsymbol{X}}_{2}\left(e^{*}\right)=\mathbf{0}$ which is clearly not optimal; because, for example, by introducing $x^{\prime} \in[0,1]^{n}$ such that $x_{1}^{\prime}=1$ and $x_{j}^{\prime}=x_{j}^{*}=0, \forall j \in\{2, \ldots, n\}$, it is concluded that $\boldsymbol{x}^{\prime}$ is feasible to (16) and $\mathbf{1}^{T} \boldsymbol{x}^{\prime}=1>0=\mathbf{1}^{T} \boldsymbol{x}^{*}$. \textbf{(a)} The proof is obtained from Lemma 13(c).\\

\section{Numerical Example}
Example 1. Consider the problem $A \otimes x=b$, where
$$
\begin{aligned}
A & =\left[\begin{array}{llllllllll}
0.81 & 0.15 & 0.65 & 0.70 & 0.43 & 0.27 & 0.75 & 0.84 & 0.35 & 0.07 \\
0.90 & 0.57 & 0.03 & 0.98 & 0.38 & 0.67 & 0.25 & 0.25 & 0.83 & 0.05 \\
0.12 & 0.95 & 0.84 & 0.27 & 0.76 & 0.14 & 0.50 & 0.81 & 0.58 & 0.53 \\
0.19 & 0.40 & 0.93 & 0.40 & 0.29 & 0.16 & 0.99 & 0.94 & 0.54 & 0.87 \\
0.63 & 0.80 & 0.67 & 0.09 & 0.45 & 0.11 & 0.89 & 0.92 & 0.91 & 0.93 \\
0.09 & 0.14 & 0.75 & 0.82 & 0.48 & 0.79 & 0.95 & 0.85 & 0.28 & 0.12 \\
0.87 & 0.42 & 0.74 & 0.69 & 0.44 & 0.95 & 0.54 & 0.19 & 0.75 & 0.56 \\
0.74 & 0.91 & 0.39 & 0.71 & 0.64 & 0.34 & 0.73 & 0.25 & 0.75 & 0.46 \\
0.95 & 0.79 & 0.65 & 0.95 & 0.70 & 0.58 & 0.14 & 0.61 & 0.38 & 0.01 \\
0.96 & 0.95 & 0.17 & 0.03 & 0.75 & 0.22 & 0.25 & 0.47 & 0.56 & 0.03
\end{array}\right] \\
b^{T} & =\left[\begin{array}{llll}
0.66,0.57,0.14,0.40,0.45,0.79,0.55,0.62,0.04,0.53
\end{array}\right] \\
c^{T} & =\left[\begin{array}{lll}
-8.36,0.92,4.11,2.36,-9.66,-8.87,5.75,-4.78,8.10,5.84
\end{array}\right]
\end{aligned}
$$
By Definition 2, we have 
$$
\begin{array}{ll}
    J_{1}=\{1,4,7,8\}, J_{2}=\{1,2,4,6,9\} \\
    J_{3}=\{2,3,4,5,6,7,8,9,10\}, J_{4}=\{2,3,4,7,8,9,10\} \\
    J_{5}=\{1,2,3,5,7,8,9,10\}, J_{6}=\{4,6,7,8\}\\
    J_{7}=\{1,3,4,6,9,10\}, J_{8}=\{1,2,4,5,7,9\}\\
J_{9}=\{1,2,3,4,5,6,7,8,9\} ~~\text{and} ~~ J_{10}=\{1,2,5,9\}
\end{array}
$$
Also, from Definition 3, we obtain \(I_{1}=\{1,3,9\}, I_{2}=\left\{i_{1}^{\prime}, i_{2}^{\prime}, i_{3}^{\prime}, i_{4}^{\prime}\right\}=\{2,4,5,6\}, I_{3}=\left\{i_{1}^{\prime \prime}, i_{2}^{\prime \prime}, i_{3}^{\prime \prime}\right\}=\{7,8,10\}\).\\

For each $i \in I_{1}$, the minimum solution $\underline{\boldsymbol{X}}(i)$ and maximum solution $\overline{\boldsymbol{X}}(i)$ of $S\left(\boldsymbol{a}_{i}, b_{i}\right)$ are obtained by Definition 4 as follows:

$$
\begin{aligned}
& \bar{X}(1)=[0.66,1,1,1,1,1,1,1,1,1] \\
& \underline{X}(1)=[0.66,0,0,0,0,0,0,0,0,0] \\
& \bar{X}(3)=[1,1,0.14,1,1,1,1,1,1,1] \\
& \underline{X}(3)=[0,0,0.14,0,0,0,0,0,0,0] \\
& \bar{X}(9)=[1,1,1,1,1,1,1,1,0.04,1] \\
& \underline{X}(9)=[0,0,0,0,0,0,0,0,0.04,0]
\end{aligned}
$$

According to Definition 5, for each $i_{k}^{\prime} \in I_{2}(k=1, \ldots, 4)$, the minimum solution $\underline{\boldsymbol{X}}(i)$ and the maximal solutions $\overline{\boldsymbol{X}}(i, 1)$ and $\overline{\boldsymbol{X}}(i, 2)$ of $S\left(\boldsymbol{a}_{i_{k}^{\prime}}, b_{i_{k}}\right)$ are attained as follows:

$$
\begin{aligned}
& \bar{X}(2,1)=[1,0.57,1,1,1,1,1,1,1,1] \\
& \overline{\boldsymbol{X}}(2,2)=[0.57,1,1,0.57,1,0.57,1,1,0.57,1] \\
& \underline{\boldsymbol{X}}(2)=[0,0.57,0,0,0,0,0,0,0,0] \\
& \overline{\boldsymbol{X}}(4,1)=[1,1,1,0.40,1,1,1,1,1,1] \\
& \overline{\boldsymbol{X}}(4,2)=[1,1,0.40,1,1,1,0.40,0.40,0.40,0.40] \\
& \underline{X}(4)=[0,0,0,0.40,0,0,0,0,0,0] \\
& \overline{\boldsymbol{X}}(5,1)=[1,1,1,1,0.45,1,1,1,1,1] \\
& \overline{\boldsymbol{X}}(5,2)=[0.45,0.45,0.45,1,1,1,0.45,0.45,0.45,0.45] \\
& \underline{\boldsymbol{X}}(5)=[0,0,0,0,0.45,0,0,0,0,0] \\
& \overline{\boldsymbol{X}}(6,1)=[1,1,1,1,1,0.79,1,1,1,1] \\
& \overline{\boldsymbol{X}}(6,2)=[1,1,1,0.79,1,1,0.79,0.79,1,1] \\
& \underline{X}(6)=[0,0,0,0,0,0.79,0,0,0,0]
\end{aligned}
$$

Also, for each $i_{k}^{\prime \prime} \in I_{3}(k=1,2,3)$, the minimal solutions $\underline{X}(i, j)$ and the maximal solutions $\overline{\boldsymbol{X}}(i, 1)$ and $\overline{\boldsymbol{X}}(i, 2)$ of $S\left(\boldsymbol{a}_{i_{k}}, b_{i_{k}}\right)$ are attained by Definition 6 as follows:

$$
\begin{aligned}
& \bar{X}(7,1)=[1,1,1,1,1,1,0.55,1,1,1] \\
& \bar{X}(7,2)=[0.55,1,0.55,0.55,1,0.55,1,1,0.55,0.55] \\
& \underline{X}(7,1)=[0.55,0,0,0,0,0,0.55,0,0,0] \\
& \underline{X}(7,3)=[0,0,0.55,0,0,0,0.55,0,0,0] \\
& \underline{X}(7,4)=[0,0,0,0.55,0,0,0.55,0,0,0] \\
& \underline{X}(7,6)=[0,0,0,0,0,0.55,0.55,0,0,0] \\
& \underline{X}(7,9)=[0,0,0,0,0,0,0.55,0,0.55,0] \\
& \underline{X}(7,10)=[0,0,0,0,0,0,0.55,0,0,0.55]
\end{aligned}
$$

$$
\begin{aligned}
& \bar{X}(8,1)=[1,1,1,1,1,1,1,0.62,1,1] \\
& \bar{X}(8,2)=[0.62,0.62,1,0.62,0.62,1,0.62,1,0.62,1] \\
& \underline{X}(8,1)=[0.62,0,0,0,0,0,0,0.62,0,0] \\
& \underline{X}(8,2)=[0,0.62,0,0,0,0,0,0.62,0,0] \\
& \underline{X}(8,4)=[0,0,0,0.62,0,0,0,0.62,0,0] \\
& \underline{X}(8,5)=[0,0,0,0,0.62,0,0,0.62,0,0] \\
& \underline{X}(8,7)=[0,0,0,0,0,0,0.62,0.62,0,0] \\
& \underline{X}(8,9)=[0,0,0,0,0,0,0,0.62,0.62,0] \\
& \bar{X}(10,1)=[1,1,1,1,1,1,1,1,1,0.53] \\
& \bar{X}(10,2)=[0.53,0.53,1,1,0.53,1,1,1,0.53,1] \\
& \underline{X}(10,1)=[0.53,0,0,0,0,0,0,0,0,0.53] \\
& \underline{X}(10,2)=[0,0.53,0,0,0,0,0,0,0,0.53] \\
& \underline{X}(10,5)=[0,0,0,0,0.53,0,0,0,0,0.53] \\
& \underline{X}(10,9)=[0,0,0,0,0,0,0,0,0.53,0.53]
\end{aligned}
$$
Therefore, from Definition 9, we have

$$
\begin{aligned}
& \overline{\boldsymbol{X}}_{1}=\min _{i \in I_{1}}\{\overline{\boldsymbol{X}}(i)\}=[0.66,1,0.14,1,1,1,1,1,0.04,1] \\
& \underline{\boldsymbol{X}}_{1}=\max _{i \in I_{1}}\{\underline{\boldsymbol{X}}(i)\}=[0.66,0,0.14,0,0,0,0,0,0.04,0] \\
& \underline{\boldsymbol{X}}_{2}=\max _{i \in I_{2}}\{\underline{\boldsymbol{X}}(i)\}=[0,0.57,0,0.40,0.45,0.79,0,0,0,0]\\
& \max \left\{\underline{\boldsymbol{X}}_{1}, \underline{\boldsymbol{X}}_{2}\right\}=[0.66,0.57,0.14,0.40,0.45,0.79,0,0,0.04,0]
\end{aligned}
$$
Based on Remark 1, we have $\left|\bar{E}_{2}\right|=2^{\left|I_{2}\right|}=2^{4}$. For example, by considering $\bar{e} \in \bar{E}_{2}$ such that $\vec{e}(2)=2, \vec{e}(4)=1, \vec{e}(5)=2$ and $\vec{e}(6)=1$, the corresponding maximal solution is obtained as follows:
$$
\begin{aligned}
\overline{\boldsymbol{X}}_{2}\left(\vec{e}^{\prime}\right) & =\min _{i \in I_{2}}\left\{\overline{\boldsymbol{X}}\left(i, \vec{e}^{\prime}(i)\right)\right\} \\
& =\min \left\{\overline{\boldsymbol{X}}(2, \vec{e}(2)), \overline{\boldsymbol{X}}(4, \vec{e}(4)), \overline{\boldsymbol{X}}(5, \vec{e}(5)), \overline{\boldsymbol{X}}\left(6, \vec{e}^{\prime}(6)\right)\right\} \\
& =\min \{\overline{\boldsymbol{X}}(2,2), \overline{\boldsymbol{X}}(4,1), \overline{\boldsymbol{X}}(5,2), \overline{\boldsymbol{X}}(6,1)\} \\
& =[0.45,0.45,0.45,0.40,1,0.57,0.45,0.45,0.45,0.45]
\end{aligned}
$$
Similarly, $\left|\bar{E}_{3}\right|=2^{\left|I_{3}\right|}=2^{3}$. For instance, if $\bar{e}^{\prime \prime} \in \bar{E}_{3}$ such that $\vec{e}^{\prime \prime}(7)=1, \vec{e}^{\prime \prime}(8)=1$ and $\vec{e}^{\prime \prime}(10)=2$ , then the corresponding maximal solution is resulted as follows:
$$
\begin{aligned}
\overline{\boldsymbol{X}}_{3}\left(\vec{e}^{\prime \prime}\right) & =\min _{i \in I_{3}}\left\{\overline{\boldsymbol{X}}\left(i, \bar{e}^{\prime \prime}(i)\right)\right\} \\
& =\min \left\{\overline{\boldsymbol{X}}\left(7, \bar{e}^{\prime \prime}(7)\right), \overline{\boldsymbol{X}}\left(8, \bar{e}^{\prime \prime}(8)\right), \overline{\boldsymbol{X}}\left(10, \bar{e}^{\prime \prime}(10)\right)\right\} \\
& =\min \{\overline{\boldsymbol{X}}(7,1), \overline{\boldsymbol{X}}(8,1), \overline{\boldsymbol{X}}(10,2)\} \\
& =[0.53,0.53,1,1,0.53,1,0.55,0.62,0.53,1]
\end{aligned}
$$
The cardinality of $\underline{E}$ is equal to $|\underline{E}|=\prod_{i \in I_{3}}\left|J_{i}\right|=\left|J_{7}\right| \times\left|J_{8}\right| \times\left|J_{10}\right|=6 \times 6 \times 4=144$. As an example, by selecting $\underline{e} \in \underline{E}$ such that $\underline{e}(7)=9, \underline{e}(8)=3$ and $\underline{e}(10)=5$, we have
$$
\begin{aligned}
\underline{\boldsymbol{X}}_{3}(\underline{e}) & =\max _{i \in I_{3}}\{\underline{\boldsymbol{X}}(i, \underline{e}(i))\} \\
& =\max \{\underline{\boldsymbol{X}}(7, \underline{e}(7)), \underline{\boldsymbol{X}}(8, \underline{e}(8)), \underline{\boldsymbol{X}}(10, \underline{e}(10))\} \\
& =\max \{\underline{\boldsymbol{X}}(7,9), \underline{\boldsymbol{X}}(8,2), \underline{\boldsymbol{X}}(10,5)\} \\
& =[0,0.62,0,0,0.53,0,0.55,0.62,0.55,0.53]
\end{aligned}
$$
Moreover, the number of all the triples $(\underline{e}, \vec{e}, \vec{e}) \in \underline{E} \times \bar{E}_{2} \times \bar{E}_{3}$ is equal to $|\underline{E}| \times\left|\bar{E}_{2}\right| \times\left|\bar{E}_{3}\right|=18432$. According to Definition 11, matrices $\overline{\boldsymbol{M}}^{1}, \overline{\boldsymbol{M}}^{2}, \overline{\boldsymbol{N}}^{1}, \overline{\boldsymbol{N}}^{2}$ and $\underline{N}$ are computed as follows:

$$
\begin{aligned}
& \overline{\boldsymbol{M}}^{1}=\left[\begin{array}{cccccccccc}
1 & 0.57 & 1 & 1 & 1 & 1 & 1 & 1 & 1 & 1 \\
1 & 1 & 1 & 0.40 & 1 & 1 & 1 & 1 & 1 & 1 \\
1 & 1 & 1 & 1 & 0.45 & 1 & 1 & 1 & 1 & 1 \\
1 & 1 & 1 & 1 & 1 & 0.79 & 1 & 1 & 1 & 1
\end{array}\right] \\
& \overline{\boldsymbol{M}}^{2}=\left[\begin{array}{ccccccccccc}
0.57 & 1 & 1 & 0.57 & 1 & 0.57 & 1 & 1 & 0.57 & 1 \\
1 & 1 & 0.40 & 1 & 1 & 1 & 0.40 & 0.40 & 0.40 & 0.40 \\
0.45 & 0.45 & 0.45 & 1 & 1 & 1 & 0.45 & 0.45 & 0.45 & 0.45 \\
1 & 1 & 1 & 0.79 & 1 & 1 & 0.79 & 0.79 & 1 & 1
\end{array}\right] \\
& \overline{\boldsymbol{N}}^{1}=\left[\begin{array}{ccccccccccccc}
1 & 1 & 1 & 1 & 1 & 1 & 0.55 & 1 & 1 & 1 \\
1 & 1 & 1 & 1 & 1 & 1 & 1 & 0.62 & 1 & 1 \\
1 & 1 & 1 & 1 & 1 & 1 & 1 & 1 & 1 & 0.53
\end{array}\right] \\
& \overline{\boldsymbol{N}}^{2}=\left[\begin{array}{cccccccccccc}
0.55 & 1 & 0.55 & 0.55 & 1 & 0.55 & 1 & 1 & 0.55 & 0.55 \\
0.62 & 0.62 & 1 & 0.62 & 0.62 & 1 & 0.62 & 1 & 0.62 & 1 \\
0.53 & 0.53 & 1 & 1 & 0.53 & 1 & 1 & 1 & 0.53 & 1
\end{array}\right] \\
& \underline{\boldsymbol{N}}=\left[\begin{array}{ccccccccccc}
0.55 & -\infty & 0.55 & 0.55 & -\infty & 0.55 & 0.55 & -\infty & 0.55 & 0.55 \\
0.62 & 0.62 & -\infty & 0.62 & 0.62 & -\infty & 0.62 & 0.62 & 0.62 & -\infty \\
0.53 & 0.53 & -\infty & -\infty & 0.53 & -\infty & -\infty & -\infty & 0.53 & 0.53
\end{array}\right]
\end{aligned}
$$
Since $\max \left\{\left(\underline{\boldsymbol{X}}_{1}\right)_{1},\left(\underline{\boldsymbol{X}}_{2}\right)_{1}\right\}=0.66>0.57=\bar{m}_{11}^{2}$ and $\max \left\{\left(\underline{\boldsymbol{X}}_{1}\right)_{1},\left(\underline{\boldsymbol{X}}_{2}\right)_{1}\right\}=0.66>0.45=\bar{m}_{31}^{2}$, by Rule 1 (Lemma 5) and Remark 3, we set $\overline{\boldsymbol{M}}_{1}^{2}=\overline{\boldsymbol{M}}_{3}^{2}=[\infty, \infty, \ldots, \infty]_{1 \times n}, \operatorname{dom}\left(\vec{e}\left(i_{1}^{\prime}\right)\right)=\operatorname{dom}\left(\vec{e}^{\prime}(2)\right)=\{1\}$ and $\operatorname{dom}\left(\vec{e}^{\prime}\left(i_{3}^{\prime}\right)\right)=\operatorname{dom}\left(\vec{e}^{\prime}(5)\right)=\{1\}$. After applying Rule $1,\left|\bar{E}_{2}\right|$ is decreased from $2^{4}=16$ to $2^{2}=4$ and the matrix $\overline{\boldsymbol{M}}^{2}$ is reduced as follows:
$$
\overline{\boldsymbol{M}}^{2}=\left[\begin{array}{cccccccccc}
\infty & \infty & \infty & \infty & \infty & \infty & \infty & \infty & \infty & \infty \\
1 & 1 & 0.40 & 1 & 1 & 1 & 0.40 & 0.40 & 0.40 & 0.40 \\
\infty & \infty & \infty & \infty & \infty & \infty & \infty & \infty & \infty & \infty \\
1 & 1 & 1 & 0.79 & 1 & 1 & 0.79 & 0.79 & 1 & 1
\end{array}\right]
$$
Similarly, by applying Rule 2 (Lemma 6), we set 
\begin{align*}                      \overline{\boldsymbol{N}}_{1}^{2}=\overline{\boldsymbol{N}}_{2}^{2}=\overline{\mathbf{N}}_{3}^{2}=[\infty, \infty, \ldots, \infty]_{1 \times n},
\end{align*}
$\operatorname{dom}\left(\vec{e}^{\prime \prime}\left(i_{1}^{\prime \prime}\right)\right)=\operatorname{dom}\left(\vec{e}^{\prime \prime}(7)\right)=\{1\}, \operatorname{dom}\left(\vec{e}^{\prime \prime}\left(i_{2}^{\prime \prime}\right)\right)=\operatorname{dom}\left(\vec{e}^{\prime \prime}(8)\right)=\{1\}$ and 
\begin{align*}                      \operatorname{dom}\left(\vec{e}^{\prime \prime}\left(i_{3}^{\prime \prime}\right)\right)=\operatorname{dom}\left(\vec{e}^{\prime \prime}(10)\right)=\{1\} ; 
\end{align*}
because, $\bar{n}_{11}^{2}=0.55, \bar{n}_{21}^{2}=0.62$ and $\bar{n}_{31}^{2}=0.53$, and therefore 
\begin{align*}                      
\bar{n}_{11}^{2}, \bar{n}_{21}^{2}, \bar{n}_{31}^{2}<\max \left\{\left(\underline{\boldsymbol{X}}_{1}\right)_{1},\left(\underline{\boldsymbol{X}}_{2}\right)_{1}\right\}.
\end{align*}
So, $\left|\bar{E}_{3}\right|$ is decreased from $2^{3}=8$ to 1 and the reduced matrix $\bar{N}^{2}$ is obtained as follows:

$$
\bar{N}^{2}=\left[\begin{array}{llllllllll}
\infty & \infty & \infty & \infty & \infty & \infty & \infty & \infty & \infty & \infty \\
\infty & \infty & \infty & \infty & \infty & \infty & \infty & \infty & \infty & \infty \\
\infty & \infty & \infty & \infty & \infty & \infty & \infty & \infty & \infty & \infty
\end{array}\right]
$$
By considering the matrix $\underline{N}$, we note that $\underline{n}_{13}=0.55$ that is greater than $\left(\overline{\boldsymbol{X}}_{1}\right)_{3}=0.14$ . Also, $\underline{n}_{19}=0.55, \underline{n}_{29}=0.62$ and $\underline{n}_{39}=0.53$ are greater than $\left(\overline{\boldsymbol{X}}_{1}\right)_{9}=0.04$. So, by applying Rule 3 (Lemma 7), we set $\underline{n}_{13}=\underline{n}_{19}=\underline{n}_{29}=\underline{n}_{39}=-\infty$ and
$$
\begin{aligned}
& \operatorname{dom}\left(\underline{e}\left(i_{1}^{\prime \prime}\right)\right)=\operatorname{dom}(\underline{e}(7))=J_{7}-\{3,9\}=\{1,4,6,10\} \\
& \operatorname{dom}\left(\underline{e}\left(i_{2}^{\prime \prime}\right)\right)=\operatorname{dom}(\underline{e}(8))=J_{8}-\{9\}=\{1,2,4,5,7\} \\
& \operatorname{dom}\left(\underline{e}\left(i_{3}^{\prime \prime}\right)\right)=\operatorname{dom}(\underline{e}(10))=J_{10}-\{9\}=\{1,2,5\}
\end{aligned}
$$
By this simplification rule, $|\underline{E}|$ is decreased from 144 to $\left|J_{7}\right| \times\left|J_{8}\right| \times\left|J_{10}\right|=4 \times 5 \times 3=60$ and the new matrix $\bar{N}^{2}$ is obtained as follows:
$$
\underline{\boldsymbol{N}}=\left[\begin{array}{cccccccccc}
0.55 & -\infty & -\infty & 0.55 & -\infty & 0.55 & 0.55 & -\infty & -\infty & 0.55 \\
0.62 & 0.62 & -\infty & 0.62 & 0.62 & -\infty & 0.62 & 0.62 & -\infty & -\infty \\
0.53 & 0.53 & -\infty & -\infty & 0.53 & -\infty & -\infty & -\infty & -\infty & 0.53
\end{array}\right]
$$
Now, consider $i_{2}^{\prime} \in I_{2} \quad\left(i_{2}^{\prime}=4\right)$ and $i_{1}^{\prime \prime} \in I_{3}\left(i_{1}^{\prime \prime}=7\right)$. It is clear that $a_{i^{\prime} i_{1}^{\prime \prime}}=a_{47}=0.99>0.4=b_{i_{2}^{\prime}}=b_{4}$ and $b_{i_{2}^{\prime}}=b_{4}=0.4<0.55=b_{i_{1}^{\prime \prime}}=b_{7}$. Hence, by applying Rule 4, we set $\overline{\boldsymbol{M}}_{2}^{2}=[\infty, \infty, \ldots, \infty]_{1 \times n}$ and $\operatorname{dom}\left(\vec{e}^{\prime}\left(i_{2}^{\prime}\right)\right)=\operatorname{dom}\left(\vec{e}^{\prime}(4)\right)=\{1\}$. After applying Rule 4, $\left|\vec{E}_{2}\right|$ is decreased from 4 to 2 and the new matrix $\overline{\boldsymbol{M}}^{2}$ is obtained as follows:
$$
\overline{\boldsymbol{M}}^{2}=\left[\begin{array}{cccccccccc}
\infty & \infty & \infty & \infty & \infty & \infty & \infty & \infty & \infty & \infty \\
\infty & \infty & \infty & \infty & \infty & \infty & \infty & \infty & \infty & \infty \\
\infty & \infty & \infty & \infty & \infty & \infty & \infty & \infty & \infty & \infty \\
1 & 1 & 1 & 0.79 & 1 & 1 & 0.79 & 0.79 & 1 & 1
\end{array}\right]
$$
In this example, Rule 5 cannot reduce the matrix $\overline{\boldsymbol{N}}^{2}$. However, by considering $i_{1}^{\prime} \in I_{2} \quad\left(i_{1}^{\prime}=2\right)$ and $i_{2}^{\prime \prime} \in I_{3}\left(i_{2}^{\prime \prime}=8\right)$, we note that $\operatorname{dom}\left(\vec{e}^{\prime}\left(i_{1}^{\prime}\right)\right)=\operatorname{dom}(\vec{e}(2))=\{1\}$, $i_{1}^{\prime}=2 \in J_{i_{2}^{\prime \prime}}=J_{8}=\{1,2,4,5,7\}$ and $b_{i_{1}^{\prime}}=b_{2}=0.57<0.62=b_{i_{2}^{\prime \prime}}=b_{8}$. Therefore, according to Rule 6, we set $\underline{n}_{2 i_{1}^{\prime}}=\underline{n}_{22}=-\infty$ and $\operatorname{dom}\left(\underline{e}\left(i_{2}^{\prime \prime}\right)\right)=\operatorname{dom}(\underline{e}(8))=J_{i_{2}^{\prime \prime}}-\left\{i_{1}^{\prime}\right\}=J_{8}-\{2\}=\{1,4,5,7\}$ . By the same argument, we set $\underline{n}_{22}=\underline{n}_{14}=\underline{n}_{24}=\underline{n}_{25}=\underline{n}_{35}=-\infty$ and
$$
\begin{aligned}
& \operatorname{dom}\left(\underline{e}\left(i_{2}^{\prime \prime}\right)\right)=\operatorname{dom}(\underline{e}(8))=J_{i_{2}^{\prime \prime}}-\left\{i_{2}^{\prime}, i_{3}^{\prime}\right\}=J_{8}-\{4,5\}=\{1,7\} \\
& \operatorname{dom}\left(\underline{e}\left(i_{1}^{\prime \prime}\right)\right)=\operatorname{dom}(\underline{e}(7))=J_{i_{1}^{\prime \prime}}-\left\{i_{2}^{\prime}\right\}=J_{7}-\{4\}=\{1,6,10\} \\
& \operatorname{dom}\left(\underline{e}\left(i_{3}^{\prime \prime}\right)\right)=\operatorname{dom}(\underline{e}(10))=J_{i_{3}^{\prime \prime}}-\left\{i_{3}^{\prime}\right\}=J_{10}-\{5\}=\{1,2\}
\end{aligned}
$$
Hence, through use of Rule $6,|\underline{E}|$ is decreased from 60 to $\left|J_{7}\right| \times\left|J_{8}\right| \times\left|J_{10}\right|=3 \times 2 \times 2=12$ and the new matrix $\bar{N}^{2}$ is reduced further as follows:
$$
\underline{\boldsymbol{N}}=\left[\begin{array}{cccccccccc}
0.55 & -\infty & -\infty & -\infty & -\infty & 0.55 & 0.55 & -\infty & -\infty & 0.55 \\
0.62 & -\infty & -\infty & -\infty & -\infty & -\infty & 0.62 & 0.62 & -\infty & -\infty \\
0.53 & 0.53 & -\infty & -\infty & -\infty & -\infty & -\infty & -\infty & -\infty & 0.53
\end{array}\right]
$$
Finally, for $i_{1}^{\prime \prime}, i_{2}^{\prime \prime} \in I_{3} \quad\left(i_{1}^{\prime \prime}=7 \quad\right.$ and $\left.i_{2}^{\prime \prime}=8\right)$ we have $\operatorname{dom}\left(\vec{e}^{\prime \prime}\left(i_{1}^{\prime \prime}\right)\right)=\operatorname{dom}\left(\vec{e}^{\prime \prime}(7)\right)=\{1\}$, $i_{1}^{\prime \prime}=7 \in J_{i_{2}^{\prime \prime}}=J_{8}=\{1,7\}$ and $b_{i_{1}^{\prime \prime}}=b_{7}=0.55<0.62=b_{i_{2}^{\prime \prime}}=b_{8}$. Thus, based on Rule 7, we set $\underline{n}_{2 i_{1}^{\prime \prime}}=\underline{n}_{27}=-\infty$ and $\operatorname{dom}\left(\underline{e}\left(i_{2}^{\prime \prime}\right)\right)=\operatorname{dom}(\underline{e}(8))=J_{i_{2}^{\prime \prime}}-\left\{i_{1}^{\prime \prime}\right\}=J_{8}-\{7\}=\{1\}$. By the same argument, we set $\underline{n}_{1,10}=-\infty$ and $\operatorname{dom}\left(\underline{e}\left(i_{1}^{\prime \prime}\right)\right)=\operatorname{dom}(\underline{e}(7))=J_{i_{1}^{\prime \prime}}\left\{i_{3}^{\prime \prime}\right\}=J_{7}-\{10\}=\{1,6\}$. So, Rule 7 decreases $|\underline{E}|$ from 12 to $\left|J_{7}\right| \times\left|J_{8}\right| \times\left|J_{10}\right|=2 \times 1 \times 2=4$ and the new matrix $\bar{N}^{2}$ is obtained as follows:\\

$\underline{\boldsymbol{N}}=\left[\begin{array}{cccccccccc}0.55 & -\infty & -\infty & -\infty & -\infty & 0.55 & 0.55 & -\infty & -\infty & -\infty \\ 0.62 & -\infty & -\infty & -\infty & -\infty & -\infty & -\infty & 0.62 & -\infty & -\infty \\ 0.53 & 0.53 & -\infty & -\infty & -\infty & -\infty & -\infty & -\infty & -\infty & 0.53\end{array}\right]$\\

Consequently, after applying the simplification rules, it follows that $\left|\bar{E}_{2}\right|=2,\left|\bar{E}_{3}\right|=1$ and $|\underline{E}|=4$. Additionally, we note that the simplification rules decreased the number of all the triples $\left(\underline{e}, \vec{e}^{\prime}, \vec{e}^{\prime \prime}\right) \in \underline{E} \times \bar{E}_{2} \times \bar{E}_{3}$ from 18432 to $|\underline{E}| \times\left|\bar{E}_{2}\right| \times\left|\bar{E}_{3}\right|=4 \times 2 \times 1=8$.

Therefore, according to Theorem 2, the feasible region $S(\boldsymbol{A}, \boldsymbol{b})$ can be found using two maximal solutions $\overline{\boldsymbol{X}}_{2}\left(\vec{e}_{1}\right), \overline{\boldsymbol{X}}_{2}\left(\vec{e}_{2}^{\prime}\right)\left(\vec{e}_{1}^{\prime}, \vec{e}_{2}^{\prime} \in \bar{E}_{2}\right)$, one maximal solution $\overline{\boldsymbol{X}}_{3}\left(\vec{e}^{\prime \prime}\right)($ $\left.\vec{e}^{\prime \prime} \in \bar{E}_{3}\right)$, and four minimal solutions $\underline{\boldsymbol{X}}_{3}\left(\underline{e}_{1}\right), \ldots, \underline{\boldsymbol{X}}_{3}\left(\underline{e}_{4}\right)\left(\underline{e}_{1}, \ldots, \underline{e}_{4} \in \underline{E}\right)$. 

These solutions are summarized as follows:

$$
\begin{aligned}
& \vec{e}_{1}^{\prime}(2)=\vec{e}_{1}(4)=\vec{e}_{1}^{\prime}(5)=\vec{e}_{1}(6)=1 \quad \Rightarrow \overline{\boldsymbol{X}}_{2}\left(\vec{e}_{1}^{\prime}\right)=[1,0.57,1,0.40,0.45,0.79,1,1,1,1] \\
& \vec{e}_{2}(2)=\vec{e}_{2}(4)=\vec{e}_{2}(5)=1, \vec{e}_{2}(6)=2 \Rightarrow \overline{\boldsymbol{X}}_{2}\left(\vec{e}_{2}^{\prime}\right)=[1,0.57,1,0.40,0.45,1,0.79,0.79,1,1] \\
& \vec{e}^{\prime \prime}(7)=\vec{e}^{\prime \prime}(8)=\vec{e}^{\prime \prime}(10)=1 ~~~~~~~~~~~ \Rightarrow \overline{\boldsymbol{X}}_{3}\left(\vec{e}^{\prime \prime}\right)=[1,1,1,1,1,1,0.55,0.62,1,0.53] \\
& \underline{e}_{1}=[1,1,1] \Rightarrow \underline{\boldsymbol{X}}_{3}\left(\underline{e}_{1}\right)=[0.62,0,0,0,0,0,0.55,0.62,0,0.53] \\
& \underline{e}_{2}=[1,1,2] \Rightarrow \underline{\boldsymbol{X}}_{3}\left(\underline{e}_{2}\right)=[0.62,0.53,0,0,0,0,0.55,0.62,0,0.53] \\
& \underline{e}_{3}=[6,1,1] \Rightarrow \underline{\boldsymbol{X}}_{3}\left(\underline{e}_{3}\right)=[0.62,0,0,0,0,0.55,0.55,0.62,0,0.53] \\
& \underline{e}_{4}=[6,1,2] \Rightarrow \underline{\boldsymbol{X}}_{3}\left(\underline{e}_{4}\right)=[0.62,0.53,0,0,0,0.55,0.55,0.62,0,0.53]
\end{aligned}
$$
So, we have
\begin{align*}
   \min \left\{\overline{\boldsymbol{X}}_{1}, \overline{\boldsymbol{X}}_{2}\left(\vec{e}_{1}^{\prime}\right), \overline{\boldsymbol{X}}_{3}\left(\vec{e}^{\prime \prime}\right)\right\}=[0.66,0.57,0.14,0.40,0.45,0.79,0.55,0.62,0.04,0.53]  
\end{align*}
\begin{align*}
    \min \left\{\overline{\boldsymbol{X}}_{1}, \overline{\boldsymbol{X}}_{2}\left(\vec{e}_{2}^{\prime}\right), \overline{\boldsymbol{X}}_{3}\left(\vec{e}^{\prime \prime}\right)\right\}=[0.66,0.57,0.14,0.40,0.45,1,0.55,0.62,0.04,0.53]
\end{align*}
Since $\min \left\{\overline{\boldsymbol{X}}_{1}, \overline{\boldsymbol{X}}_{2}\left(\vec{e}_{1}^{\prime}\right), \overline{\boldsymbol{X}}_{3}\left(\vec{e}^{\prime \prime}\right)\right\} \leq \min \left\{\overline{\boldsymbol{X}}_{1}, \overline{\boldsymbol{X}}_{2}\left(\vec{e}_{2}^{\prime}\right), \overline{\boldsymbol{X}}_{3}\left(\vec{e}^{\prime \prime}\right)\right\}$, then 
\begin{align*}
    \min \left\{\overline{\boldsymbol{X}}_{1}, \overline{\boldsymbol{X}}_{2}\left(\vec{e}_{2}^{\prime}\right), \overline{\boldsymbol{X}}_{3}\left(\vec{e}^{\prime \prime}\right)\right\}
\end{align*}
is the unique maximum solution of $S(\boldsymbol{A}, \boldsymbol{b})$. Moreover, since $\underline{\boldsymbol{X}}_{3}\left(\underline{e}_{1}\right) \leq \underline{\boldsymbol{X}}_{3}\left(\underline{e}_{k}\right), k=2,3,4$ , then 
\begin{align*}
    \max \left\{\underline{\boldsymbol{X}}_{1}, \underline{\boldsymbol{X}}_{2}, \underline{\boldsymbol{X}}_{3}\left(\underline{e}_{1}\right)\right\}=[0.66,0.57,0.14,0.40,0.45,0.79,0.55,0.62,0.04,0.53]
\end{align*}
is the unique minimal (minimum) solution of $S(\boldsymbol{A}, \boldsymbol{b})$. Hence, according to Theorem 2, we have $\quad \boldsymbol{S}(\boldsymbol{A}, \boldsymbol{b})=\left[\max \left\{\underline{\boldsymbol{X}}_{1}, \underline{\boldsymbol{X}}_{2}, \underline{\boldsymbol{X}}_{3}\left(\underline{e}_{1}\right)\right\}, \min \left\{\overline{\boldsymbol{X}}_{1}, \overline{\boldsymbol{X}}_{2}\left(\vec{e}_{2}^{\prime}\right), \overline{\boldsymbol{X}}_{3}\left(\vec{e}^{\prime \prime}\right)\right\}\right] . \quad$ Additionally, $J^{+}=\{2,3,4,7,9,10\}$ and $J^{-}=\{1,5,6,8\}$. Thus, from Theorem 4, the optimal solution of the problem is obtained as 
\begin{align*}
    \boldsymbol{x}_{e^{*}}=[0.66,0.57,0.14,0.40,0.45,1,0.55,0.62,0.04,0.53]
\end{align*}
where $e^{*}=\left(\underline{e}_{1}, \vec{e}_{2}, \vec{e}^{\prime \prime}\right)$, and then $\boldsymbol{c}^{T} \boldsymbol{x}_{e^{*}}=-13.07$.

\section{Conclusion}
In this paper, an algorithm was proposed for finding a global optimal solution of linear objective problems subjected to a novel system of fuzzy relation equations defined with minimum t-norm. It was proved that the feasible solutions set is completely resolved by a finite number of closed convex cells. Some necessary and sufficient conditions were also presented to determine the feasibility of the problem. In addition, seven simplification rules were proposed to speed up solution finding. It was proved that the problem has a binary optimum if all the coefficients and the values of the variables are binary. Moreover, it was shown that the well-known minimum vertex cover problem is actually a special case of the problem presented in this paper.

\section{Acknowledgment}
We are very grateful to the anonymous referees for their comments and suggestions, which were very helpful in improving the paper.

\section{Reference}
[1]. C. W. Chang, B. S. Shieh, Linear optimization problem constrained by fuzzy max-min relation equations, Information Sciences 234 (2013) 71-79.

[2]. L. Chen, P. P. Wang, Fuzzy relation equations (i): the general and specialized solving algorithms, Soft Computing 6 (5) (2002) 428-435.

[3]. L. Chen, P. P. Wang, Fuzzy relation equations (ii): the branch-point-solutions and the categorized minimal solutions, Soft Computing 11 (1) (2007) 33-40.

[4]. M. Cornejo, D. Lobo, J. Medina, Bipolar fuzzy relation equations based on product t-norm, in: Proceedings of 2017 IEEE International Conference on Fuzzy Systems, 2017.

[5]. S. Dempe, A. Ruziyeva, On the calculation of a membership function for the solution of a fuzzy linear optimization problem, Fuzzy Sets and Systems 188 (2012) 58-67.

[6]. A. Di Nola, S. Sessa, W. Pedrycz, E. Sanchez, Fuzzy relational Equations and their applications in knowledge engineering, Dordrecht: Kluwer Academic Press, 1989.

[7]. D. Dubey, S. Chandra, A. Mehra, Fuzzy linear programming under interval uncertainty based on IFS representation, Fuzzy Sets and Systems 188 (2012) 68-87. [8]. D. Dubois, H. Prade, An introduction to bipolar representations of information and preference, International Journal of Intelligent Systems 23(8) (2008) 866-877.

[9]. D. Dubois, H. Prade, An overview of the asymmetric bipolar representation of positive and negative information in possibility theory, Fuzzy Sets and Systems 160 (2009) 1355-1366.

[10]. Y. R. Fan, G. H. Huang, A. L. Yang, Generalized fuzzy linear programming for decision making under uncertainty: Feasibility of fuzzy solutions and solving approach, Information Sciences 241 (2013) 12-27.

[11]. S.C. Fang, G. Li, Solving fuzzy relational equations with a linear objective function, Fuzzy Sets and Systems 103 (1999) 107-113.

[12]. S. Freson, B. De Baets, H. De Meyer, Linear optimization with bipolar max-min constraints, Information Sciences 234 (2013) 3-15.

[13]. A. Ghodousian, E. Khorram, Fuzzy linear optimization in the presence of the fuzzy relation inequality constraints with max-min composition, Information Sciences 178 (2008) 501-519.

[14]. A. Ghodousian, E. Khorram, Linear optimization with an arbitrary fuzzy relational inequality, Fuzzy Sets and Systems 206 (2012) 89-102.

[15]. A. Ghodousian, M. Raeisian Parvari, A modified PSO algorithm for linear optimization problem subject to the generalized fuzzy relational inequalities with fuzzy constraints (FRI-FC), Information Sciences 418-419 (2017) 317-345.

[16]. A. Ghodousian, A. Babalhavaeji, An efficient genetic algorithm for solving nonlinear optimization problems defined with fuzzy relational equations and max-Lukasiewicz composition, Applied Soft Computing 69 (2018) 475-492.

[17]. A. Ghodousian, M. Naeeimib, A. Babalhavaeji, Nonlinear optimization problem subjected to fuzzy relational equations defined by Dubois-Prade family of t-norms, Computers \& Industrial Engineering 119 (2018) 167-180.

[18]. A. Ghodousian, Optimization of linear problems subjected to the intersection of two fuzzy relational inequalities defined by Dubois-Prade family of t-norms, Information Sciences 503 (2019) 291-306.

[19]. A. Ghodousian, F. Samie Yousefi, Linear optimization problem subjected to fuzzy relational equations and fuzzy constraints, Iranian Journal of Fuzzy Systems ? (2022) ???-???.

[20]. F. F. Guo, L. P. Pang, D. Meng, Z. Q. Xia, An algorithm for solving optimization problems with fuzzy relational inequality constraints, Information Sciences 252 ( 2013) 20-31.

[21]. S. M. Guu, Y. K. Wu, Minimizing a linear objective function under a max-t-norm fuzzy relational equation constraint, Fuzzy Sets and Systems 161 (2010) 285-297. [22]. S. M. Guu, Y. K. Wu, Minimizing a linear objective function with fuzzy relation equation constraints, Fuzzy Optimization and Decision Making 12 (2002) 1568-4539.

[23]. H. C. Lee, S. M. Guu, On the optimal three-tier multimedia streaming services, Fuzzy Optimization and Decision Making 2(1) (2002) 31-39.

[24]. P. K. Li, S. C. Fang, On the resolution and optimization of a system of fuzzy relational equations with sup-t composition, Fuzzy Optimization and Decision Making 7 (2008) 169-214.

[25]. J. X. Li, S. J. Yang, Fuzzy relation inequalities about the data transmission mechanism in bittorrent-like peer-to-peer file sharing systems, in: Proceedings of the $9^{\text {th }}$ International Conference on Fuzzy Systems and Knowledge discovery (FSKD 2012), pp. 452-456.

[26]. P. Li, Q. Jin, On the resolution of bipolar max-min equations, Kybernetika 52(4) (2016) 514530.

[27]. P. Li, Y. Liu, Linear optimization with bipolar fuzzy relational equation constraints using lukasiewicz triangular norm, Soft Computing 18 (2014) 1399-1404.

[28]. J. L. Lin, Y. K. Wu, S. M. Guu, On fuzzy relational equations and the covering problem, Information Sciences 181 (2011) 2951-2963.

[29]. J. L. Lin, On the relation between fuzzy max-archimedean t-norm relational equations and the covering problem, Fuzzy Sets and Systems 160 (2009) 2328-2344.

[30]. C. C. Liu, Y. Y. Lur, Y. K. Wu, Linear optimization of bipolar fuzzy relational equations with max-Łukasiewicz composition, Information Sciences 360 (2016) 149-162.

[31]. J. LU, S.C. Fang, Solving nonlinear optimization problems with fuzzy relation equations consraints, Fuzzy Sets and Systems 119 (2001) 1-20.

[32]. J. Loetamonphong, S. C. Fang, Optimization of fuzzy relation equations with max-product composition, Fuzzy Sets and Systems 118 (2001) 509-517.

[33]. A. V. Markovskii, On the relation between equations with max-product composition and the covering problem, Fuzzy Sets and Systems 153 (2005) 261-273.

[34]. M. Mizumoto, H. J. Zimmermann, Comparison of fuzzy reasoning method, Fuzzy Sets and Systems 8 (1982) 253-283.

[35]. W. Pedrycz, Granular Computing: Analysis and Design of Intelligent Systems, CRC Press, Boca Raton, 2013.

[36]. W. Pedrycz, Fuzzy relational equations with generalized connectives and their applications, Fuzzy Sets and Systems 10 (1983) 185-201.

[37]. W. Pedrycz, Solving fuzzy relational equations through logical filtering, Fuzzy Sets and Systems 81 (1996) 355-363. [38]. X. B. Qu, X. P. Wang, Man-hua. H. Lei, Conditions under which the solution sets of fuzzy relational equations over complete Brouwerian lattices form lattices, Fuzzy Sets and Systems 234 (2014) 34-45.

[39]. X. B. Qu, X. P. Wang, Minimization of linear objective functions under the constraints expressed by a system of fuzzy relation equations, Information Sciences 178 (2008) 3482-3490.

[40]. E. Sanchez, Solution in composite fuzzy relation equations: application to medical diagnosis in Brouwerian logic, in: M.M. Gupta. G.N. Saridis, B.R. Games (Eds.), Fuzzy Automata and Decision Processes, North-Holland, New York, 1977, pp. 221-234.

[41]. E. Sanchez, Resolution of eigen fuzzy sets equations, Fuzzy Sets and Systems 1 (1978) 69-75.

[42]. B. S. Shieh, Infinite fuzzy relation equations with continuous t-norms, Information Sciences 178 (2008) 1961-1967.

[43]. B. S. Shieh, Minimizing a linear objective function under a fuzzy max-t-norm relation equation constraint, Information Sciences 181 (2011) 832-841.

[44]. G.B. Stamou, S.G. Tzafestas, Resolution of composite fuzzy relation equations based on Archimedean triangular norms, Fuzzy Sets and Systems 120 (2001) 395-407.

[45]. F. Sun, Conditions for the existence of the least solution and minimal solutions to fuzzy relation equations over complete Brouwerian lattices, Information Sciences 205 (2012) 86-92.

[46]. F. Sun, X. P. Wang, x. B. Qu, Minimal join decompositions and their applications to fuzzy relation equations over complete Brouwerian lattices, Information Sciences 224 (2013) 143-151.

[47]. Y. K. Wu, Optimization of fuzzy relational equations with max-av composition, Information Sciences 177 (2007) 4216-4229.

[48]. Y. K. Wu, S. M. Guu, Minimizing a linear function under a fuzzy max-min relational equation constraints, Fuzzy Sets and Systems 150 (2005) 147-162.

[49]. Y. K. Wu, S. M. Guu, An efficient procedure for solving a fuzzy relation equation with maxArchimedean t-norm composition, IEEE Transactions on Fuzzy Systems 16 (2008) 73-84.

[50]. Y. K. Wu, S. M. Guu, J. Y. Liu, Reducing the search space of a linear fractional programming problem under fuzzy relational equations with max-Archimedean t-norm composition, Fuzzy Sets and Systems 159 (2008) 3347-3359.

[51]. Q. Q. Xiong, X. P. Wang, Fuzzy relational equations on complete Brouwerian lattices, Information Sciences 193 (2012) 141-152.

[52]. X.P. Yang, Resolution of bipolar fuzzy relation equations with max-Lukasiewicz composition, Fuzzy Sets and Systems (2020). [53]. X. P. Yang, X. G. Zhou, B. Y. Cao, Latticized linear programming subject to max-product fuzzy relation inequalities with application in wireless communication, Information Sciences 358-359 (2016) 44-55.

[54]. S. J. Yang, An algorithm for minimizing a linear objective function subject to the fuzzy relation inequalities with addition-min composition, Fuzzy Sets and Systems 255 (2014) 41-51.

[55]. J. Zhou, Y. Yu, Y. Liu, Y. Zhang, Solving nonlinear optimization problems with bipolar fuzzy relational equations constraints, Journal of Inequalities and Applications 126 (2016) 1-10.

\end{large}
\end{document}